\documentclass[12pt,english]{amsproc}
\usepackage[T1]{fontenc}
\usepackage[utf8]{inputenc}
\usepackage{mathrsfs, amsmath, amstext, amsthm, amssymb}
\usepackage[all]{xy}
\usepackage{babel}
\usepackage{caption}
\makeatletter

\numberwithin{equation}{section}
\numberwithin{figure}{section}

\theoremstyle{plain}
\newtheorem{thm}{\protect\theoremname}[section]
\theoremstyle{remark}
\newtheorem{rem}[thm]{\protect\remarkname}
\theoremstyle{plain}
\newtheorem{prop}[thm]{\protect\propositionname}
\theoremstyle{definition}
\newtheorem{defn}[thm]{\protect\definitionname}
\theoremstyle{remark}
\newtheorem{example}[thm]{\protect\examplename}
\theoremstyle{plain}
\newtheorem{lem}[thm]{\protect\lemmaname}
\theoremstyle{plain}
\newtheorem{cor}[thm]{\protect\corollaryname}

\usepackage{graphicx}
\usepackage{wrapfig}
\usepackage{fancyhdr}
\usepackage{extarrows}
\usepackage{geometry}
\newtheoremstyle{notestyle}
 {\topsep}
 {\topsep}
 {\normalfont}
 {}
 {\bfseries}
 {.}
 {.5em}
 {}

\theoremstyle{notestyle}

\makeatother

\providecommand{\definitionname}{Definition}
\providecommand{\lemmaname}{Lemma}
\providecommand{\propositionname}{Proposition}
\providecommand{\remarkname}{Remark}
\providecommand{\theoremname}{Theorem}
\providecommand{\examplename}{Example}
\providecommand{\corollaryname}{Corollary}

\usepackage{hyperref}

\title{On the Euclidean Distance Degree of Shallow Polynomial Networks}
\author{Giacomo Graziani}
\address{%
Department of Civil, Environmental and Architectural Engineering\\
University of Padua\\
Via F. Marzolo 9\\
35131 Padova, Italy}

\email{giacomo.graziani@unipd.it}



\makeatletter
\newcommand{\unlistedsection}[1]{%
  \begingroup
  \let\@tocwrite\@gobbletwo
  \section*{#1}%
  \endgroup
}
\newcommand{\unlistedsubsection}[1]{%
  \begingroup
  \let\@tocwrite\@gobbletwo
  \subsection*{#1}%
  \endgroup
}
\makeatother
\begin{document}

\global\long\def\AA{\mathbb{A}}%
\global\long\def\PP{\mathbb{P}}%
\global\long\def\CC{\mathbb{C}}%
\global\long\def\RR{\mathbb{R}}%
\global\long\def\QQ{\mathbb{Q}}%

\global\long\def\XX{\mathfrak{X}}%

\global\long\def\R{\mathrm{R}}%

\global\long\def\V{\boldsymbol{V}}%
\global\long\def\r{\boldsymbol{\underline{r}}}%
\global\long\def\a{\boldsymbol{\underline{\alpha}}}%

\global\long\def\kk{\Bbbk}%

\global\long\def\p{^{\prime}}%
\global\long\def\d{^{\vee}}%
\global\long\def\ot{\otimes}%
\global\long\def\op{\oplus}%
\global\long\def\inv{^{\ast}}%
\global\long\def\su{\subseteq}%
\global\long\def\i{_{\mathrm{in}}}%
\global\long\def\o{_{\mathrm{out}}}%
\global\long\def\e{\emptyset}%

\global\long\def\F{\mathcal{F}}%
\global\long\def\I{\mathcal{I}}%
\global\long\def\X{\mathcal{X}}%
\global\long\def\M{\mathcal{M}}%
\global\long\def\N{\mathcal{N}}%
\global\long\def\W{\mathcal{W}}%
\global\long\def\G{\mathcal{G}}%
\global\long\def\OO{\mathcal{O}}%
\global\long\def\UU{\mathcal{U}}%
\global\long\def\EE{\mathcal{E}}%
\global\long\def\C{\mathcal{C}}%
\global\long\def\Q{\mathcal{Q}}%
\global\long\def\P{\mathcal{P}}%
\global\long\def\T{\mathcal{T}}%
\global\long\def\Y{\mathcal{Y}}
\global\long\def\c{^{\circ}}%
\global\long\def\K{\mathcal{K}}%
\global\long\def\H{\mathcal{H}}%
\global\long\def\D{\mathcal{D}}%

\global\long\def\cl{\mathrm{cl}}%

\global\long\def\ga{\alpha}%
\global\long\def\gb{\beta}%
\global\long\def\gs{\sigma}%
\global\long\def\gl{\lambda}%
\newcommand{\gt}{\theta}

\global\long\def\g{\mathfrak{g}}%
\global\long\def\Om{\mathfrak{W}}%
\global\long\def\WW{\mathfrak{W}}%

\global\long\def\Sym{\mathrm{Sym}}%
\global\long\def\Hom{\mathrm{Hom}}%
\global\long\def\Seg{\mathrm{Seg}}%
\global\long\def\GL{\mathrm{GL}}%
\global\long\def\ED{\mathrm{ED}}
\global\long\def\ged{\mathrm{gED}}%
\global\long\def\adim{\mathrm{adim}}%
\global\long\def\edim{\mathrm{edim}}%
\global\long\def\Aut{\mathrm{Aut}}%
\global\long\def\Net{\mathrm{Net}}%
\global\long\def\En{\mathrm{En}}%
\global\long\def\Gr{\mathrm{Gr}}%
\global\long\def\pr{\mathrm{pr}}%
\global\long\def\Sub{\mathcal{D}_{\mathrm{neur}}}%
\global\long\def\Da{\mathcal{D}_{\mathrm{arch}}}%
\global\long\def\Of{\Om^{\mathrm{full}}}%
\global\long\def\Lie{\mathrm{Lie}}%
\global\long\def\rk{\mathrm{rk}}%
\global\long\def\Nash{\mathrm{Nash}}%
\global\long\def\Hilb{\mathrm{Hilb}}%
\global\long\def\Td{\mathrm{Td}}%
\global\long\def\td{\mathrm{td}}%
\global\long\def\ch{\mathrm{ch}}%

\begin{abstract}
We study the projective geometry of shallow polynomial neural networks, namely families of homogeneous polynomial maps realised by a single hidden layer with monomial activation. The closure of the set of maps realised by such a network is a projective variety, its neurovariety. A natural measure of the algebraic complexity of the corresponding best-approximation problem is the generic Euclidean distance degree $\mathrm{gED}$, which counts the complex critical points of the squared distance from a general point.

For fixed input dimension, width and activation degree, we prove that $\mathrm{gED}$ is eventually polynomial in the output dimension. Its degree and leading coefficient are governed by an auxiliary variety parametrising the spaces spanned by powers of linear forms. When the activation degree is at least the width, this extends: on the same stable range, $\mathrm{gED}$ is a polynomial in the output dimension and the activation degree jointly.

Our proof constructs a proper birational model of the conormal variety over a base that does not depend on the output dimension, and expresses $\mathrm{gED}$ through Chern classes on it. We give closed formulas in two basic cases and determine the leading term for width two.
\end{abstract}

\maketitle
\tableofcontents 
\newpage

\section{Introduction}

The algebraic geometry of neurovarieties has attracted increasing attention in recent years. This is due to the central role played by algebraic and semi-algebraic models in deep learning, supported by density and universal approximation results \cite{SM17,SM18}, by the existence of global invariants describing the expressive power of an architecture \cite{KTB19,KLW24}, by the relative tractability of the singularities these varieties exhibit - more delicate to handle in a differential-geometric setting \cite{Wat09} - and by the appearance of programmatic works such as \cite{MSMTK25}.

Two natural invariants in this setting are the dimension of the neurovariety and the algebraic degree of the associated optimisation problem. The former provides a geometric measure of expressivity, while the latter controls the number of complex critical points arising in a general nearest-point problem.  This algebraic degree is the so called Euclidean distance degree: given a scalar product $q$ on the space of functions, $\ED_q(X)$ counts the complex critical points of the squared distance from a general target. It therefore measures the algebraic complexity of the corresponding nearest-point problem and bounds the number of isolated real critical points. It depends sensitively on $q$, and its value may decrease for scalar products compatible with the symmetries of $X$ \cite{KMQS25}; for a general scalar product, however, it attains a fixed maximal value, a purely projective invariant of $X$ given by the sum of its polar degrees, the generic Euclidean distance degree $\ged\left(X\right)$ \cite{DHOST}.

Computing $\ged$ for a neurovariety is hard, and explicit values are currently available mainly for individual architectures, including the quadratic network $\CC^{2}\to\CC^{2}\to\CC^{d_2}$ studied in \cite[Theorem 6.5]{KLW24} and several rank-one varieties considered in \cite[Table 2]{KMQS25}. This suggests studying its behaviour along families in which one or more architectural parameters vary. In this paper we fix a shallow polynomial network of width $k$ and activation degree $r$ with input of dimension $n$, and we let the output dimension $m$ grow: we prove that $\ged$ is eventually a polynomial in $m$, whose degree and leading coefficient are governed by an auxiliary variety attached to the network, and that for $r\ge k$ it is jointly a polynomial in $m$ and $r$.

The main obstacle to computing $\ged\left(X_{k,r}^{n,m}\right)$ is that the invariant is defined by the conormal cycle $\left[\N_{X_m}\right]$, which lives in $\PP\left(\M_m\right)\times\PP\left(\M_m\right)\d$, an ambient space whose dimension grows with $m$. Our first idea is to find a birational model $\C_m$ of $\N_{X_m}$, constructed over a base $B$ independent of $m$, and of $r$ as well whenever $r\ge k$, and read all the dependence on the parameters inside characteristic classes on $B$. Everything rests on the auxiliary variety $\Sigma=\Sigma_{k,r}^{n}\su\Gr_k\left(S_r\right)$ of $k$-dimensional spaces of $r$-th powers. It does not depend on $m$, by construction; and for $r\ge k$ it does not depend on $r$ either, since Proposition \ref{prop:BaseFissa} identifies it with $\Hilb_{k}^{\mathrm{sm}}\PP\left(V\d\right)$, whose dimension is $k\left(n-1\right)$ by Proposition \ref{prop:DimensioneSigma}.

\subsection*{Main results}%~\\
Our two main results are the following. The first fixes the architecture  $\left(n,k,r\right)$ and lets the output dimension $m$ grow: by Theorem \ref{thm:PolinomialitaStabileInM} there is a polynomial $P_{n,k,r}$ of degree exactly $b=\dim\left(\Sigma\right)$ with
\[
\ged\left(X_{k,r}^{n,m}\right)=P_{n,k,r}\left(m\right)\qquad\text{for every }m\ge kn ,
\]
whose leading coefficient is $4^{b}\deg\Sigma/b!$. The second lets the activation degree $r$ vary as well: by Theorem \ref{thm:PolinomialitaCongiunta}, for $r\ge k$ there is a polynomial $Q_{n,k}\in\QQ\left[m,r\right]$, of degree at most $b$ in each variable separately, with $\ged\left(X_{k,r}^{n,m}\right)=Q_{n,k}\left(m,r\right)$ on the whole quadrant $r\ge k$, $m\ge kn$. The theorems are actually proved with the sharper threshold $m\ge k+b$, in which $k+b\le kn$ because $b\le k\left(n-1\right)$; its optimality is discussed in Remark \ref{rem:DueCondizioni} and in Appendix \ref{sec:esempi}. The construction behind both is of independent interest, in two ways. In Proposition \ref{prop:ProprietaCm} we build a proper birational model of the conormal variety of $X_{k,r}^{n,m}$ over a base which depends neither on $m$ nor, for $r\ge k$, on $r$, hence turning the study of a family of cycles in ambient spaces of growing dimension into the intersection theory on a single fixed smooth variety. Moreover the argument never uses what $\Sigma$ is: Theorem \ref{thm:GeneraleSigma} states the first result for an arbitrary irreducible $\Sigma\su\Gr_{k}\left(S\right)$ and the variety $\bigcup_{L\in\Sigma}\PP\left(W\ot L\right)$ swept out by the corresponding subspaces. This is what makes it portable. Under a mild genericity assumption the neurovariety of a polynomial network of arbitrary depth has that shape, with $\Sigma$ the variety of the spaces of polynomials computed by the last hidden layer, so that the generic $\ED$-degree is a polynomial in the output dimension in that generality as well; what remains open is the geometry of that variety, and we discuss in Section \ref{sec:Finale}. Along the way we evaluate the resulting formula in closed  form for $k=1$ and for $n=k=2$ (Propositions \ref{prop:FormulaCasoK1} and \ref{prop:FormulaCasoBinario}), and we compute the degree of $\Sigma_{2,r}^{n}$ for every $n$ (Lemma \ref{lem:GradoSigma}), which gives the leading term of $P_{n,2,r}$. The slice $r=2$ of the second formula is the determinantal case treated in \cite[Theorem 6.5]{KLW24}, whose stated range of validity we sharpen in Remark \ref{rem:ConfrontoKLW}. Finally, Appendix \ref{app:Genericity} contains a self-contained proof, valid over any algebraically closed field of characteristic different from $2$, of the genericity statement underlying Proposition \ref{prop:ProprietaGed}$\left(1\right)$.

\subsection*{Content of the paper}%~\\
In Section \ref{sec:neurovarieta} we introduce the neurovariety $X=X_{k,r}^{n,m}$ and the bound on the width $k$ which makes its geometry non-trivial, that is, the architecture non-filling (Remark \ref{rem:kMinoreN}). We then introduce the incidence variety $\Y=\PP_{\Sigma}\left(W\ot\UU\right)\to\PP\left(\M\right)$ and use it to isolate a smooth open subset $X\c$ containing a real point (Remark \ref{rem:PuntiRealiLisci}), to relate $\dim X$ to $\dim\left(\Sigma\right)$ (Proposition \ref{prop:XcircLiscio}), and to describe the affine tangent space of $X$ along $X\c$ (Corollary \ref{cor:SpazioTangente}).\newline
In Section \ref{sec:gED} we briefly recall the definition of Euclidean distance degree of an embedded variety associated to a nondegenerate symmetric bilinear form and introduce in Definition \ref{def:DefinizioneGradiPolariConormale} a purely projective invariant defined by means of the conormal class $\left[\N_{X}\right]$ of the variety. In Proposition \ref{prop:ProprietaGed} we recall the relation between these invariants and in Proposition \ref{prop:FormulaChernLiscio} we give a formula that reduces the computation to a function of Chern classes when $X$ is smooth. This formula will allow us to perform computations in Appendix \ref{sec:esempiK1}. \newline
Section \ref{sec:conormale} is the technical heart of the paper, where we construct the birational model $\C_m$ of the conormal variety announced above. In view of the description of the tangent space given in Section \ref{sec:neurovarieta}, the condition defining $\N_{X_m}$ over $X_m\c$ can be written over $\Sigma$ (Lemma \ref{lem:CriterioConormale}). This allows us to realise $\N_{X_m}$ as the image of the zero locus $\C_m$ of a section of a vector bundle over a desingularisation $B$ of the Nash blow-up of $\Sigma$: in Proposition \ref{prop:ProprietaCm} we prove that, for $m\ge k+b$, the scheme $\C_m$ is an integral local complete intersection of codimension $b$ and that the natural morphism $\Phi:\C_m\to\N_{X_m}$ is proper and birational. Since $B$ does not depend on the parameters, the whole dependence on $m$, and later on $r$, is carried by the bundles over it, which is what makes the conormal cycles of different networks comparable.\newline
In Section \ref{sec:polinomialita} we put the previous constructions together and obtain formula \eqref{eq:gEDFormulaBuona}, which expresses $\ged\left(X_m\right)$ as an integral over $B$. We then deduce in Theorem \ref{thm:PolinomialitaStabileInM} that $\ged\left(X_m\right)$ is given by a polynomial in $m$ for $m\ge k+b$, and we compute its leading term: it is governed by the degree of $\Sigma$ in the Plücker embedding, which we make explicit for $k=2$ in Appendix \ref{sec:esempiBinario}.\newline
In Section \ref{sec:polinomialitaCongiunta} we let $r$ vary as well. We prove that $\Sigma_{k,r}^{n}$ is isomorphic to $\Hilb_{k}^{\mathrm{sm}}\PP\left(V\d\right)$, which does not involve $r$ (Proposition \ref{prop:BaseFissa}), so that for $r\ge k$ a single base $B$ serves every $r$ and the whole dependence on the activation degree is concentrated in the bundle $\EE_{r}=p_{\ast}\xi\inv\OO\left(r\right)$; the bound $r\ge k$ cannot be dispensed with (Remark \ref{rem:SogliaRK}). A singular Riemann-Roch argument then shows that the Chern classes of $\EE_{r}$ are polynomial in $r$ (Lemma \ref{lem:ChernPolinomialiInR}), and feeding this into \eqref{eq:gEDFormulaBuona} gives Theorem \ref{thm:PolinomialitaCongiunta}: for $r\ge k$ and $m\ge kn$ the generic $\ED$-degree is a polynomial in $m$ and $r$ jointly, of degree at most $b$ in each variable separately.\newline
In Section \ref{sec:Finale} we observe that none of the preceding argument used the description of $\Sigma$: Theorem \ref{thm:GeneraleSigma} holds for an arbitrary irreducible subvariety of a Grassmannian, and we indicate what this gives for polynomial networks of arbitrary depth.

\unlistedsection{Acknowledgments}
I am grateful to Kathlén Kohn and Jan Draisma for bringing important references to my attention and for their comments on the first draft of this paper. I also thank the anonymous referee of a previous version of this work for pointing out several errors and for comments that led to substantial improvements of the manuscript.

\newpage

\section{Neurovarieties}\label{sec:neurovarieta}

Throughout, $V$ and $W$ are real vector spaces of dimensions $n\ge2$ and $m$, and $k,r\ge1$ are integers. The shallow polynomial network of width $k$ and activation degree $r$ with input $V$ and output $W$ is the family of polynomial maps
\[
f_{\gt_{1},\gt_{2}}=\gt_{2}\circ\rho_r\circ\gt_{1}:V\to W,\qquad
\left(\gt_{1},\gt_{2}\right)\in\Hom\left(V,\RR^{k}\right)\times\Hom\left(\RR^{k},W\right),
\]
where $\rho_r$ raises each coordinate of $\RR^{k}$ to the $r$-th power; we consider networks  with no bias, all neurons carrying the same activation. Writing $\ell_{1},\dots,\ell_k\in V\d$ for the rows of $\gt_{1}$ and $w_{1},\dots,w_k\in W$ for the columns of $\gt_{2}$, one has $f_{\gt_{1},\gt_{2}}\left(v\right)=\sum_iw_i\,\ell_i\left(v\right)^{r}$, so that the set of maps realised by the network identifies with
\[
\F=\Big\{ \sum_{i=1}^{k}w_i\ot\ell_i^{r}\ \Big\vert\ w_i\in W,\ \ell_i\in V\d\Big\}
\su\M=W\ot S_r,\qquad S_r=\Sym^{r}\left(V\d\right).
\]
All the geometric constructions of this paper are performed after extending scalars to $\CC$: we write $V_{\CC},W_{\CC},\M_{\CC}$ for the complexifications and we call neurovariety the Zariski closure
\begin{equation}\label{eq:Neurovarieta}
X=X_{k,r}^{n,m}=\overline{\PP\left(\F\right)}
=\gs_k\left[\Seg\left(\PP\left(W_{\CC}\right)\times v_r\PP\left(V_{\CC}\d\right)\right)\right]
\su\PP\left(\M_{\CC}\right),
\end{equation}
the $k$-th secant variety of a Segre-Veronese variety; it is a projective variety defined over $\RR$, and $\F$ is contained in the affine cone over the real points of $X$. See \cite{KTB19,KLW24}. From now on we drop the subscript $\CC$, all varieties being considered over  $\CC$ unless the real structure is explicitly involved. We set
\[
N=\dim S_r=\binom{n+r-1}{r},\qquad\dim\PP\left(\M \right)=mN-1,
\]
and for $T\in\M $ we denote with $T^{\flat}:W\d\to S_r$ the associated linear map, so that $T=\sum_iw_i\ot\ell_i^{r}$ has $\mathrm{Im}\left(T^{\flat}\right)\su\left\langle \ell_{1}^{r}, \dots,\ell_k^{r}\right\rangle $.
% \begin{example}
% For $r=1$ the network is linear and $X$ is the variety of $m\times n$ matrices of rank at most $k$, while for $k=1$ the variety $X=\Seg\left(\PP\left(W\right)\times v_r\PP\left(V\d\right)\right)$ is smooth.
% \end{example}

Let
\[
\Sigma=\Sigma_{k,r}^{n}=\overline{\left\{ \left\langle \ell_{1}^{r},\dots,\ell_k^{r}\right\rangle \,\vert\,\ell_i\in V\d,\ \ell_{1}^{r},\dots,\ell_k^{r}\text{ independent}\right\} }\su\Gr_k\left(S_r\right)
\]
be the variety of $k$-dimensional spaces of $r$-th powers, an irreducible variety of dimension 
\[
\dim\left(\Sigma\right)=:b\le k\left(n-1\right).
\]

\begin{rem}\label{rem:kMinoreN}
When $m>k$, we see $X$ is a proper subvariety of $\PP\left(\M\right)$ if and only if $k<N$. Indeed the $r$-th powers span $S_r$, so if $k\ge N$ we may choose $\ell_{1},\dots,\ell_N\in V\d$ with $\ell_{1}^{r},\dots,\ell_N^{r}$ a basis of $S_r$ and write every $T\in W\ot S_r$ as $T=\sum_{i=1}^{N}w_i\ot\ell_i^{r}$ hence the architecture is filling and $X=\PP\left(\M\right)$. On the other hand if $k<N$ then $\min\left(m,N\right)\ge k+1$, so under the identification $\M\simeq\Hom\left(W\d,S_r\right)$ there is a $T$ with $\rk T^{\flat}=k+1$; but on $X$ we have $\rk T^{\flat}\le k$ so $\left[T\right]\notin X$. In order to avoid exceptions in the statements we shall always assume $k<N$.
\end{rem}

Let $\UU$ be the restriction to $\Sigma$ of the tautological subbundle of $\Gr_k\left(S_r\right)$ and let
\[
\Y=\PP_{\Sigma}\left(W\ot\UU\right)=\left\{ \left(L,\left[T\right]\right)\,\vert\,T\in W\ot L\right\} \su\Sigma\times\PP\left(\M \right),
\]
we denote the projections with $p:\Y\to\Sigma$ and $\mu:\Y\to\PP\left(\M \right)$. Then $p$ is a projective bundle with fibre $\PP\left(W\ot L\right)$, hence $\Y$ is irreducible of dimension $km+b-1$.

Our next aim is to isolate an open subset $X\c\su X$ over which $\mu$ is an isomorphism, so as to obtain the intrinsic description \eqref{eq:TangenteIntrinseco} of the affine tangent space of $X$: this is the form in which the geometry of $X$ enters the conormal construction of Section \ref{sec:conormale}.

\begin{lem}\label{lem:IncidenzaPropria}
We have that
\[
X =\bigcup_{L\in\Sigma }\PP\left(W\ot L\right).
\]
Moreover if $\rk T^{\flat}=k$ then $\mathrm{Im}\left(T^{\flat}\right)\in\Sigma $ and $\left(\mu \right)^{-1}\left(\left[T\right]\right)=\left\{ \left(\mathrm{Im}\,T^{\flat},\left[T\right]\right)\right\}$.
\end{lem}

\begin{proof}
The variety $\Y $ is closed in $\Sigma \times\PP\left(\M\right)$ and $\Sigma $ is projective, hence $\mu $ is proper and its image is closed. Let $\Sigma_{0}\su\Sigma$ denote the dense open subset of the spaces $L=\left\langle \ell_{1}^{r},\dots,\ell_k^{r}\right\rangle $ with $\ell_{1}^{r}, \dots, \ell_k^{r}$ independent. If $L\in\Sigma_{0}$, every $T\in W\ot L$ is of the form $\sum_iw_i\ot\ell_i^{r}$, so that $\mu \left(\Y_{\vert\Sigma_{0}}\right) \su X$ and hence $\mu\left(\Y\right) \su X$. Since $\mu\left(\Y\right)$ is closed and contains every tensor $\sum_{i=1}^{k}w_i\ot\ell_i^{r}$ with $\ell_{1}^{r},\dots,\ell_k^{r}$ independent and these are dense in $X$ we have  $\mu\left(\Y\right)=X$. The last statement is clear.
\end{proof}

We write $\Sigma_{\mathrm{reg}}$ for the smooth locus of $\Sigma$ and set
\[
X\c=\left\{\left[T\right]\in X\,\vert\,\rk T^{\flat}=k,\ \mathrm{Im}\left(T^{\flat}\right)\in\Sigma_{\mathrm{reg}}\right\} ,
\qquad\Y\c=\mu^{-1}\left(X\c\right)
\]

\begin{prop}\label{prop:XcircLiscio}
Let $m\ge k$, then $X \c$ is a nonempty open subset of $X$,
\[
\Y\c=\left\{ \left(L,\left[T\right]\right)\in\Y \,\vert\,L\in\Sigma_{\mathrm{reg}},\ \rk T^{\flat}=k\right\}
\]
and $\mu$ restricts to an isomorphism $\Y\c\overset{\sim}{\to}X \c$. In particular $X \c$ is smooth, irreducible and dense in $X$, so that $X \c\su X_{\mathrm{reg}}$ and $\dim X=km+b-1$.
\end{prop}

\begin{proof}
The set $\left\{ \rk T^{\flat}\ge k\right\} $ is open, so that
\[
X^{\left(k\right)}:=X\cap\left\{ \rk T^{\flat}\ge k\right\}
=\left\{ \left[T\right]\in X\,\vert\,\rk T^{\flat}=k\right\}
\]
is an open subvariety of $X$. The image of the tautological map
\[
W\d\ot\OO_{X^{\left(k\right)}}\left(-1\right)\longrightarrow S_r\ot\OO_{X^{\left(k\right)}}
\]
is a subbundle of rank $k$ and it defines a morphism $\gamma:X^{\left(k\right)}\longrightarrow\Gr_k\left(S_r\right)$ with $\gamma\left(\left[T\right]\right)=\mathrm{Im}\left(T^{\flat}\right)$ which takes values in $\Sigma $ in view of  Lemma \ref{lem:IncidenzaPropria}. Therefore  $X\c=\gamma^{-1}\left(\Sigma_{\mathrm{reg}}\right)$ is open in $X^{\left(k\right)}$, hence in $X$. It is nonempty: since $m\ge k$, for any $L\in\Sigma_{\mathrm{reg}}$ there exists a surjection $W\d\twoheadrightarrow L$, that is a tensor $T\in W\ot L$ with $\mathrm{Im}\left(T^{\flat}\right)=L$, and $\left[T\right]\in X$ in view again of  Lemma \ref{lem:IncidenzaPropria}. The description of $\Y\c$ and the equalities $\mu\circ g=\mathrm{id}$ and $g\circ\mu=\mathrm{id}$ for
\[
g=\left(\gamma,\mathrm{id}\right):X \c\to\Sigma_{\mathrm{reg}}\times\PP\left(\M\right)
\]
follow from Lemma \ref{lem:IncidenzaPropria}; note that $g$ takes values in $\Y\c$, so $\mu$ and $g$ are mutually inverse isomorphisms. Finally $\Y\c$ is a nonempty open subset of the projective bundle $\PP_{\Sigma_{\mathrm{reg}}}\left(W\ot\UU \right)$ over the smooth variety $\Sigma_{\mathrm{reg}}$, hence it is smooth, irreducible, of dimension $b+km-1$ and dense in $\Y$; everything is transported to $X \c$ by the isomorphism.
\end{proof}

\begin{rem}\label{rem:PuntiRealiLisci}
Assume $m\ge k$ and $k<N$. Then $X$ has smooth real points. Indeed, choose $\ell_{1},\dots,\ell_k\in V\d$ and $w_{1},\dots,w_k\in W$ general and real, and set $T=\sum_iw_i\ot\ell_i^{r}$. The powers $\ell_{1}^{r},\dots,\ell_k^{r}$ are then linearly independent, and so are the $w_i$ since $m\ge k$, so that $\rk T^{\flat}=k$ and $\mathrm{Im}\left(T^{\flat}\right)=\left\langle \ell_{1}^{r},\dots,\ell_k^{r}\right\rangle $. Moreover this span lies in $\Sigma_{\mathrm{reg}}$: on the open subset of $\PP\left(V\d\right)^{k}$ where the powers $\ell_{1}^{r},\dots,\ell_k^{r}$ are independent, the preimage of $\Sigma_{\mathrm{sing}}$ under $\left(\ell_{1},\dots,\ell_k\right)\mapsto\left\langle \ell_{1}^{r},\dots,\ell_k^{r}\right\rangle $ is a proper closed subset defined over $\RR$, and a proper closed subset cannot contain all its real points, these being Zariski dense. Hence $T$ is a real point of $X\c$, which is contained in $X_{\mathrm{reg}}$ by Proposition \ref{prop:XcircLiscio}. For $k\ge N$ the architecture is filling and $X=\PP\left(\M_{\CC}\right)$, so the conclusion is trivial.
\end{rem}

\begin{cor}\label{cor:ConoAffine}
Let $m\ge k$ and denote by $\mathrm{Tot}\left(W\ot\UU \right)$ the total space of the bundle $W\ot\UU $ over $\Sigma $. Then the map
\[
\left(L,T\right)\longmapsto T
\]
induces an isomorphism between the open subset $\left\{ \left(L,T\right)\,\vert\,L\in\Sigma_{\mathrm{reg}},\ \rk T^{\flat}=k\right\}$ and the punctured affine cone over $X\c$.
\end{cor}

\begin{proof}
It is the argument of Proposition \ref{prop:XcircLiscio} one level up: the map is bijective onto the punctured affine cone over $X\c$ by Lemma \ref{lem:IncidenzaPropria}, and its inverse $T\mapsto\left(\gamma\left(\left[T\right]\right),T\right)$ is a morphism, $\gamma$ being the morphism constructed there.
\end{proof}

\begin{cor}\label{cor:SpazioTangente}
Let $m\ge k$, let $T\in X\c$ and $L=\mathrm{Im}\left(T^{\flat}\right)$. Put $Q_L =S_r/L$, let $q_L :S_r\to Q_L $ be the projection and let
\[
\vartheta_T :T_L \Sigma\longrightarrow W\ot Q_L ,\qquad \vartheta_T \left(\xi\right)=\left(W\ot\xi\right)\left(T\right),
\]
where we identify $T_L \Sigma$ with a subspace of $\Hom\left(L,Q_L \right)$ and view $T\in W\ot L$. Then $\vartheta_T $ is injective,
\begin{equation}\label{eq:TangenteIntrinseco} 
\widehat{T}_T X=\left(W\ot q_L \right)^{-1}\left(\vartheta_T \left(T_L \Sigma\right)\right),
\end{equation}
and $\dim \widehat{T}_T X = mk+b$.
\end{cor}

\begin{proof}
Write $T=\sum_{a=1}^{k}w_a\ot e_a$ with $e_{1},\dots,e_k$ a basis of $L$. Then $T^{\flat}\left(\varphi\right)=\sum_a\varphi\left(w_a\right)e_a$, so that $T^{\flat}$ is onto $L$ if and only if $w_{1},\dots,w_k$ are linearly independent. Since $\rk T^{\flat}=k$, from $\vartheta_T\left(\xi\right)=\sum_aw_a\ot\xi\left(e_a\right)=0$ it follows that $\xi\left(e_a\right)=0$ for every $a$, hence $\xi=0$ and $\vartheta_T $ is injective.

We compute the Zariski tangent space of the affine cone by dual numbers $\CC\left[\varepsilon\right]$. In view of Corollary \ref{cor:ConoAffine} an element $T\p\in W\ot S_r$ is tangent to the affine cone over $X$ at $T$ if and only if there is a $\CC\left[\varepsilon\right]$-point $L_{\varepsilon}$ of $\Sigma$ over $L$ such that $T+\varepsilon T\p\in W\ot L_{\varepsilon}$. Such an $L_{\varepsilon}$ is the same thing as an element $\xi\in T_L \Sigma\su\Hom\left(L,Q_L \right)$, realised as the $\CC\left[\varepsilon\right]$-submodule of $S_r\ot\CC\left[\varepsilon\right]$ generated by the $\ell+\varepsilon\tilde{\xi}\left(\ell\right)$, $\ell\in L$, for any lift $\tilde{\xi}\in\Hom\left(L,S_r\right)$ of $\xi$; the submodule does not depend on the lift, since changing $\tilde{\xi}$ by a map $L\to L$ amounts to composing with an automorphism of $L\ot\CC\left[\varepsilon\right]$. A general element of $W\ot L_{\varepsilon}$ is then
\[
\sum_a\left(A_a+\varepsilon B_a\right)\ot\left(e_a+\varepsilon\tilde{\xi}\left(e_a\right)\right) = \sum_aA_a\ot e_a+\varepsilon\Big(\sum_aB_a\ot e_a+\sum_aA_a\ot\tilde{\xi}\left(e_a\right)\Big)
\]
with $A_a,B_a\in W$ and imposing that its $\varepsilon^{0}$-part be $T$ gives $A_a=w_a$, so that $T\p=\sum_aB_a\ot e_a+\left(W\ot\tilde{\xi}\right)\left(T\right)$ with the $B_a\in W$ arbitrary. Therefore
\[
\widehat{T}_T X=W\ot L+\left\{ \left(W\ot\tilde{\xi}\right)\left(T\right)\,\vert\,
\xi\in T_L \Sigma\right\} ,
\]
the right-hand side being independent of the choice of the lifts. Applying $W\ot q_L $, whose kernel is $W\ot L$, this set is exactly the preimage of $\vartheta_T \left(T_L \Sigma\right)$, which is \eqref{eq:TangenteIntrinseco}. Finally $W\ot q_L $ is surjective with kernel of dimension $mk$, and $\vartheta_T \left(T_L \Sigma\right)$ has dimension $b$ because $\vartheta_T $ is injective and $L\in\Sigma_{\mathrm{reg}}$, whence $\dim \widehat{T}_T X=mk+b$.
\end{proof}

We close this section with the computation of the dimension of $\Sigma$, which will play a central role in what follows. Since we are over $\CC $, the assignment $f\ot\ell^{r}\mapsto f\left(\ell\right)$ induces a well defined pairing $S_{r}\d\otimes S_{r}\to\CC$ that identifies $H^{0}\left(\PP\left(V\d\right),\OO\left(r\right)\right)=\Sym^{r}V$ with $S_{r}\d $. For a subscheme $Z\su\PP\left(V\d\right)$ we denote by $I_{Z}$ its saturated ideal.

\begin{lem}\label{lem:Regolarita}
Let $Z\su\PP\left(V\d\right)$ have length $s$ and let $r\ge s-1$. Then the restriction $S_{r}\d\to H^{0}\left(\OO_{Z}\left(r\right)\right)$ is surjective; equivalently $\dim\left(I_{Z}\right)_{r}=N-s$, and $\dim\left\langle v_{r}\left(Z\right)\right\rangle =s-1$.
\end{lem}

\begin{proof}
The Hilbert polynomial of $Z$ is the constant $s$, whose Macaulay decomposition $$s=\binom t0+\dots+\binom{t-s+1}0$$ has $s$ terms; hence the Gotzmann number of $Z$ is $s$ and $I_{Z}$ is $s$-regular \cite[Theorem 4.3.2]{BH}, \cite[Proposition 2.1.2]{BGL13}. The assertion is then \cite[Lemma 2.1.3]{BGL13}.
\end{proof}

\begin{lem}\label{lem:Apolarita}
Let $r\ge k$ and let $\ell_{1},\dots,\ell_k\in V\d$ be pairwise non proportional. Then $\ell_{1}^{r},\dots,\ell_k^{r}$ are linearly independent and, setting $L=\left\langle \ell_{1}^{r},\dots,\ell_k^{r}\right\rangle \su S_r$,
\[
\left\{ \left[\ell\right]\in\PP\left(V\d\right)\,\vert\,\ell^{r}\in L\right\} = \left\{ \left[ \ell_{1} \right],\dots,\left[ \ell_k \right]\right\} .
\]
\end{lem}

\begin{proof}
Let $Z=\left\{ \left[ \ell_{1} \right],\dots,\left[ \ell_k \right]\right\} $, a reduced subscheme of length $k$. By Lemma \ref{lem:Regolarita} applied to $Z$, with $r\ge k-1$, one has $\dim\left\langle v_r\left(Z\right)\right\rangle =k-1$, that is $\ell_{1}^{r},\dots,\ell_k^{r}$ are linearly independent. If $\left[\ell\right]\notin Z$, then $Z\cup\left\{ \left[\ell\right]\right\} $ has length $k+1$ and Lemma \ref{lem:Regolarita} applied to it, with $r\ge k$, gives that $\ell_{1}^{r},\dots,\ell_k^{r},\ell^{r}$ are linearly independent; in particular $\ell^{r}\notin L$.
\end{proof}

\begin{prop}\label{prop:DimensioneSigma}
Let $r\ge k$. Then $b=\dim\left(\Sigma\right) =k\left(n-1\right)$.
\end{prop}

\begin{proof}
By definition $\Sigma $ is the closure of the image of the rational map
\[
\varphi:\PP\left(V\d\right)^{k}\dashrightarrow\Gr_k\left(S_r\right),\qquad
\varphi\left(\ell_{1},\dots,\ell_k\right)=\left\langle \ell_{1}^{r},\dots,\ell_k^{r}\right\rangle ,
\]
which is defined wherever the powers are independent, a nonempty open subset since $k\le N$ and the $r$-th powers span $S_r$. Hence $b\le k\left(n-1\right)$. In view of Lemma \ref{lem:Apolarita} the fibre of $\varphi$ over a general point of the image is contained in $\left\{ \left[ \ell_{1} \right],\dots,\left[ \ell_k \right]\right\} ^{k}$, hence finite, so $\varphi$ is generically finite onto its image and $b=k\left(n-1\right)$.
\end{proof}

The hypothesis $r\ge k$ is used only in the second application of Lemma \ref{lem:Regolarita}, to a scheme of length $k+1$, and it is optimal: for $r=k-1$ any $k$ points of $\PP^{1}$ span the whole of $\Sym^{k-1}\CC^{2}$, so that $L=S_{k-1}$ and every $\ell^{k-1}$ lies in $L$.

\begin{rem}\label{rem:DimensioneAttesaEkMinoreN}
It follows from Propositions \ref{prop:XcircLiscio} and \ref{prop:DimensioneSigma} that $\dim X=k\left(m+n-1\right)-1$ for $m,r\ge k$; this was well known, compare \cite{BCC11}. In particular for $k<N$ $\dim X$ is the expected dimension
\[
\min\left\{ k\left(m+n-1\right)-1,\ m\binom{n+r-1}{r}-1\right\} .
\]
The hypothesis $m>k$ of Remark \ref{rem:kMinoreN} cannot be dropped: for $\left(n,k,r\right) = \left(2,2,2\right)$ and $m=2$ one has $k=2<3=N$, and yet $X_{2,2}^{2,2}=\PP^{5}$, see Remark \ref{rem:ConfrontoKLW}.
\end{rem}

\section{The generic $\ED$-degree}\label{sec:gED}

Let $E$ be a complex vector space and let
\[
\mathcal{S}\left(E\right)
=\left\{q\in\Sym^{2}E\d\,\vert\,\det q\neq0\right\}
\]
denote the space of nondegenerate symmetric bilinear forms on $E$. Let $q\in\mathcal{S}\left(E\right)$, let $Z\su E$ be a reduced and irreducible variety not contained in the isotropic quadric $q\left( z,z\right)=0$ of $q$, and let $u\in E$ be general. The function
\[
z\longmapsto q\left(z-u,z-u\right)
\]
has finitely many critical points on $Z_{\mathrm{reg}}$, and their number is the $\ED$-degree $\ED_q\left(Z\right)$ of $Z$; see \cite[\S2]{KMQS25}. For a projective variety $X\su\PP\left(E\right)$ not contained in the isotropic quadric of $q$ we set
\[
\ED_q\left(X\right)=\ED_q\left(\widehat{X}\right),
\]
where $\widehat{X}\su E$ is the affine cone over $X$; see \cite[Lemma 2.8 and Corollary 2.9]{DHOST}.

Since Hermitian forms play no role for us, given $q\in\mathcal{S}\left(E\right)$, we say by an abuse of terminology that it is a scalar product if it is the complex bilinear extension of a positive definite symmetric bilinear form. If $q$ is a scalar product and $X$ is defined over $\RR$ with a smooth real point, then, since the isotropic quadric of $q$ has no real points, $X$ is clearly not contained in the isotropic quadric of $q$ and $\ED_q\left(X\right)$ is the number of critical points of the squared distance from a general $u$; see \cite[Theorem 4.1]{DHOST}. In our setting $X$ is a neurovariety and has a smooth real point by Remark \ref{rem:PuntiRealiLisci}, while $u$ is a function to be learned. Thus $\ED_q(X)$ measures the algebraic complexity of the nearest-point problem on $\widehat X$. When training is formulated as the minimisation of a quadratic loss in coefficient space, it gives an upper bound for the number of isolated critical points for a general target; compare \cite[\S6]{KLW24}.

In this paper we will mainly deal with an upper bound $\ged\left(X\right)$ that can be given to the family $\ED_q\left(X\right)$ for varying $q$. The interest in it lies in the fact that it can be expressed purely in terms of projective invariants of $X$. From now on $X\subsetneq\PP^M$ is a reduced and irreducible variety of dimension $d$; the assumption $X\neq\PP^M$ is needed, since for $X=\PP^M$ no hyperplane contains a tangent space and the conormal variety below is empty, while $\ED_q\left(X\right)=1$ for every $q$.

\begin{defn}\label{def:DefinizioneGradiPolariConormale}
The conormal variety of $X$ is
\[
\N_X =\overline{\left\{ \left(x,H\right)\in\PP^M\times\left(\PP^M\right)\d\,\vert\,x\in X_{\mathrm{reg}},\ T_x X\su H\right\} },
\]
an irreducible variety of dimension $M-1$, the fibre over a smooth point $x$ being the space of hyperplanes containing $T_x X$. Its image in $\left(\PP^{M}\right)\d$ is the dual variety $X\d$. Denoting by $x,y$ the hyperplane classes of the two factors, the polar degrees $\delta_i\left(X\right)$ are defined by the multidegree of the conormal cycle,
\begin{equation}\label{eq:DefinizioneGradiPolariConormale}
\left[\N_X \right]=\sum_{i=0}^{d}\delta_i\left(X\right)x^{M-i}y^{i+1}\in
A^{\bullet}\left(\PP^{M}\times\left(\PP^{M}\right)\d\right),
\end{equation}
with the convention $\delta_i\left(X\right)=0$ for $i>d$, and the generic Euclidean distance degree of $X$ is
\begin{equation}\label{eq:DefinizioneGEDGradiPolari}
\ged\left(X\right)=\sum_{i=0}^{d}\delta_i\left(X\right).
\end{equation}
\end{defn}

The name is justified by the following statement, which also explains in which sense $\ged$ controls $\ED_q$ when defined; we refer to \cite{Piene78} for polar classes and to \cite{DHOST,KMQS25} for the $\ED$-degree.

\begin{prop}\label{prop:ProprietaGed}
Let $X\subsetneq\PP\left(E\right)$ be reduced and irreducible. Then
\begin{enumerate}
\item $\ED_q\left(X\right)=\ged\left(X\right)=\ged\left(X\d\right)$ for every $q\in\mathcal{S}\left(E\right)$ such that $\N_X\cap\Delta_q=\e$, where $\Delta_q$ is the diagonal induced by $q$, and such $q$ form a nonempty Zariski open subset of $\mathcal{S}\left(E\right)$;
\item $\ED_q\left(X\right)\le\ged\left(X\right)$ for every $q\in\mathcal{S}\left(E\right)$ whose isotropic quadric does not contain $X$, such $q$ forming an open subset of $\mathcal{S}\left(E\right)$ on which $q\mapsto\ED_q\left(X\right)$ is lower semicontinuous;
\item a general real positive definite scalar product satisfies $\ED_q\left(X\right)=\ged\left(X\right)$, the real positive cone being Zariski dense in $\Sym^{2}E\d$ and therefore meeting the open subset of $\left(1\right)$.
\end{enumerate}
\end{prop}

\begin{proof}
Part $\left(1\right)$ is \cite[Theorem 2.10]{KMQS25} together with the remark following it, see also \cite[Theorem 5.4]{DHOST} and Proposition \ref{prop:GenericitaGenerale} for the proof of openness; 
part $\left(2\right)$ is \cite[Theorem 3.11]{KMQS25}, which gives that the generic value of $\ED_q\left(X\right)$ is its maximum;
part $\left(3\right)$ follows from $\left(1\right)$, the cone of real positive definite forms being a nonempty Euclidean open subset of the real vector space $\Sym ^2 E\d_\RR$, hence Zariski dense in its complexification.
\end{proof}

In Appendix \ref{app:Genericity} we give a direct proof of the genericity statement underlying Proposition \ref{prop:ProprietaGed}$\left(1\right)$, namely that $\N_X \cap \Delta_q=\e$ for a general $q$, where $\Delta _q$ denotes the diagonal induced by $q$. The argument is purely projective and works over any algebraically closed field of characteristic different from $2$, which may be of independent interest.
 
When $X$ is smooth, $\ged\left(X\right)$ is expressed by the Chern classes of $X$; this is the form in which we shall use it in Appendix \ref{sec:esempiK1}.

\begin{prop}\label{prop:FormulaChernLiscio}
Let $X\subsetneq\PP^{M}$ be a smooth variety of dimension $d$ with hyperplane class $H$. Then
\begin{equation}\label{eq:FormulaChernLiscio}
\ged\left(X\right)=\sum_{i=0}^{d}\left(-1\right)^{i}\left(2^{d+1-i}-1\right)
\int_X c_i\left(T_X \right)\cdot H^{d-i}.
\end{equation}
\end{prop}

\begin{proof}
See \cite[Theorem 5.8]{DHOST} or \cite[Theorem 2.11]{KMQS25}. 
\end{proof}
We recall that for singular $X$ the same formula holds with the Chern classes replaced by the Chern-Mather classes of $X$, see \cite{Piene78} and \cite[Proposition 2.9]{Alu18}.

\section{A birational model of the conormal variety}\label{sec:conormale}

From now on $n,k,r$ are fixed, with (recall Remark \ref{rem:kMinoreN})
\[
k<N=\dim S_r = \binom{n+r-1}{r}
\]
and $m$ varies: write $W$ for a real vector space of dimension $m\ge k$, and
\[
\M_m=W\ot S_r,\qquad X_m\su\PP\left(\M_m\right),\qquad
\mu_m:\Y_m=\PP_{\Sigma}\left(W\ot\UU\right)\to\PP\left(\M_m\right)
\]
for the objects of Section \ref{sec:neurovarieta} attached to $W$. The point of what follows is that $\Sigma$, and with it $b$ and $N$, do not depend on $m$: the whole dependence on the output dimension is carried by $W$, of dimension $m$, through the bundles $W\ot\UU$ and their duals. We shall exploit this by realising the conormal variety of $X_m$, whose ambient space $\PP\left(\M_m\right)\times\PP\left(\M_m\right)\d$ does depend on $m$, as the image of a scheme $\C_m$ living over a base built once and for all from $\Sigma$. By Remark \ref{rem:PuntiRealiLisci} the variety $X_m$ has a smooth real point, so that the hypotheses of Section \ref{sec:gED} are satisfied.

\subsection{Tangent spaces and the conormal condition}%~\\
Our goal now is to express $\N_{X_m}$ using the vanishing of linear and bilinear objects (Lemma \ref{lem:CriterioConormale}), this will allow us in Section \ref{sec:costruzioneCm} to construct a birational model $\C_m$ of $\N_{X_m}$ as the vanishing of a global section of a vector bundle.

In view of Proposition \ref{prop:XcircLiscio} the subset $X_m\c$ is open, dense and contained in $\left(X_m\right)_{\mathrm{reg}}$, so that the conormal variety of $X_m$ may be computed over $X_m\c$.

\begin{lem}\label{lem:ConormaleSuXcirc}
Let $m\left(N-k\right)>b$, equivalently $X_m \subsetneq\PP \left(\M_m\right)$, and set 
$$\N\c_{X_m}=\left\{ \left(\left[T\right],\left[Y\right]\right)\in X_m\c\times\PP\left(\M_m\d\right)\,\vert\,Y_{\vert\widehat{T}_T X_m}=0\right\}.$$
Then $\N\c_{X_m}\to X_m\c$ is the projectivisation of a vector bundle of rank $m\left(N-k\right)-b>0$; in particular it is smooth and irreducible of dimension $mN-2$, and $\N_{X_m}=\overline{\N\c_{X_m}}$.
\end{lem}

\begin{proof}
In view of Corollary \ref{cor:SpazioTangente} the affine tangent space $\widehat{T}_T X_m\su\M_m$ has dimension $km+b$ at every $T\in X_m\c$, so that
\[
\left(\widehat{T}_T X_m\right)^{\perp}
=\left\{ Y\in\M_m\d\,\vert\,Y_{\vert\widehat{T}_T X_m}=0\right\}
\]
has constant dimension $mN-km-b>0$. These subspaces are therefore the fibres of a vector bundle of that rank on $X_m\c$, and $\N\c_{X_m}$ is its projectivisation, of dimension
$$\left(km+b-1\right)+\left(mN-km-b-1\right)=mN-2.$$ 
Being a projective bundle over an irreducible variety, the projectivised conormal bundle over $\left(X_m\right) _{\mathrm{reg}}$ is irreducible, and $\N\c_{X_m}$ is a nonempty open subset of it in view of Proposition \ref{prop:XcircLiscio}, hence dense; taking closures gives the last assertion.
\end{proof}

By Proposition \ref{prop:XcircLiscio} one has $\dim X_m=km+b-1$, so that $m\left(N-k\right)\ge b$ always holds, with equality if and only if $X_m=\PP\left(\M_m\right)$, in which case $\N_{X_m}=\e$. The hypothesis is thus a condition on the pair $\left(k,m\right)$ and not on $k$ alone, and will be automatic under the assumptions of Proposition \ref{prop:ProprietaCm}, since $k<N$ and $m\ge k+b$ give $m\left(N-k\right)-b\ge m-b\ge k>0$.

Fix now $L\in\Sigma_{\mathrm{reg}}$ and put $Q_L =S_r/L$. Identifying $\M_m\d=W\d\ot S_r\d$ the functionals vanishing on $W\ot L$ are exactly those in $W\d\ot Q\d_L$, and we let
\[
\kappa_L :\left(W\ot L\right)\ot\left(W\d\ot Q\d_L \right)\longrightarrow L\ot Q\d_L \simeq\Hom\left(L,Q_L \right)\d\longrightarrow\left(T_L \Sigma\right)\d
\]
be the contraction of the factors $W$ and $W\d$ followed by the map dual to the inclusion 
$$T_L \Sigma\su\Hom\left(L,Q_L \right).$$

\begin{lem}\label{lem:CriterioConormale}
Let $\left[T\right]\in X_m\c$, let $L=\mathrm{Im}\left(T^{\flat}\right)$ and let $Y\in\M_m\d$. Then $Y$ vanishes on $\widehat{T}_T X_m$ if and only if
\[
Y\in W\d\ot Q\d_L \qquad\text{and}\qquad\kappa_L \left(T\ot Y\right)=0 .
\]
\end{lem}

\begin{proof}
In view of Corollary \ref{cor:SpazioTangente} we have $\widehat{T}_T X_m=\left(W\ot q_L \right)^{-1}\left(\vartheta_T \left(T_L \Sigma\right)\right)$, so $Y$ annihilates it if and only if it annihilates $W\ot L=\ker\left(W\ot q_L \right)$, that is $Y\in W\d\ot Q\d_L $, and the induced functional $\overline{Y}$ on $W\ot Q_L $ annihilates $\vartheta_T \left(T_L \Sigma\right)$. The latter condition reads $\overline{Y}\circ\vartheta_T =0$, and $\overline{Y}\circ\vartheta_T =\kappa_L \left(T\ot Y\right)$ by definition of $\kappa_L $: writing $T=\sum_aw_a\ot e_a$ and $Y=\sum_c\varphi_c\ot\psi_c$, both sides pair with $\xi\in T_L \Sigma$ giving $\sum_{a,c}\varphi_c\left(w_a\right)\psi_c\left(\xi\left(e_a\right)\right)$.
\end{proof}

Both conditions involve only the triple $\left(L,T,Y\right)$ with $L\in\Sigma_{\mathrm{reg}}$, $T\in W\ot L$ and $Y\in W\d\ot Q\d_L $, the second being the vanishing of an expression bilinear in $\left(T,Y\right)$ with values in $\left(T_L \Sigma\right)\d$. They therefore make sense over any variety mapping to $\Sigma$ on which the tangent spaces $T_L \Sigma$ extend to a vector bundle, which is what the Nash blow-up does.

\subsection{The variety $\C_m$}\label{sec:costruzioneCm}%~\\
We want to apply now the criterion of the previous section and construct a birational model of $\mathcal{N}_{X_m}$ over a base that does not depend on $m$. To do this we first have to move to a variety on which the tangent spaces of $\Sigma$ extend to a vector bundle. We shall in addition want the base to be smooth, so that the ambient space of $\C_m$ is smooth and $\C_m$, being cut out by as many equations as its codimension, is a local complete intersection: this is what makes it integral in Proposition \ref{prop:ProprietaCm}, and therefore its fundamental class computable in Section \ref{sec:polinomialita}. We therefore replace the Nash blow-up of $\Sigma$ by a desingularisation of it, over which we build the ambient space and the section whose zero locus is $\C_m$. Let $\Nash\left(\Sigma\right)$ be the Nash blow-up of $\Sigma$ and let
\begin{equation}\label{eq:MappaNu}
\nu:B\overset{\pi}{\longrightarrow}\Nash\left(\Sigma\right)\to\Sigma
\end{equation}
where $\pi$ is a proper resolution of singularities of $\Nash\left(\Sigma\right)$ which is an isomorphism over the inverse image of $\Sigma _\mathrm{reg}$, in particular $B$ is smooth and
$$B\c=\nu^{-1}\left(\Sigma_{\mathrm{reg}}\right)\to\Sigma_{\mathrm{reg}}$$
is an isomorphism. Over $B$ we have the pullbacks of the tautological subbundle $\EE$ and quotient $\Q$ of $\Gr_k\left(S_r\right)$, fitting in $0\to\EE\to S_r\ot\OO_B\to\Q\to0$, together with the Nash bundle $\T\su\EE\d\ot\Q$, of rank $b$, whose fibre at a point over $L \in \Sigma_{\mathrm{reg}}$ is $T_L \Sigma$. We point out that no part of this  construction depends on $m$.

Set $\P_m=\PP_B\left(W\ot\EE\right)\times_B\PP_B\left(W\d\ot\Q\d\right)$, it is a smooth variety of dimension $mN+b-2$ whose points are triples $\left(z,\left[T\right],\left[Y\right]\right)$ with $T\in W\ot\EE_z$ and $Y\in W\d\ot\Q\d_z$. Given the relative version of $\kappa_L $,
\[
\left(W\ot\EE\right)\ot\left(W\d\ot\Q\d\right)\to\EE\ot\Q\d\to\T\d ,
\]
the usual yoga of sections of vector bundles determines a global section 
\[
s_m\in H^{0}\left(\P_m,\T\d\ot\OO_{\P_m}\left(1,1\right)\right) ,
\]
and we set $\C_m=Z\left(s_m\right)\su\P_m$.

\begin{prop}
\label{prop:ProprietaCm}Let $k<\dim S_r=N$ and $m\ge k+b$. Define
\[
\P_m^{\circ}=\left\{ \left(z,\left[T\right],\left[Y\right]\right)\in\P_m\,\vert\,z\in B^{\circ},\rk T^{\flat}=k\right\} 
\]
and $\C_m^{\circ}=\C_m\cap\P_m^{\circ}$. Then
\begin{enumerate}
\item The map
\[
\begin{aligned}\Phi:\C_m & \to\PP\left(\mathcal{M}_m\right)\times\PP\left(\mathcal{M}_m\d\right)\\
\left(z,\left[T\right],\left[Y\right]\right) & \mapsto\left(\left[T\right],\left[Y\right]\right)
\end{aligned}
\]
is well-defined and gives an isomorphism $\C_m^{\circ}\simeq\mathcal{N}_{X_m}^{\circ}$;
\item $\C_m^{\circ}$ is a dense open subset of $\C_m$;
\item the image of $\Phi $ is $\mathcal{N}_{X_m}$ and the morphism $\Phi :\C_m\to\N_{X_m}$ is proper and birational.
\end{enumerate}
Moreover $\C_m$ is an integral scheme and is a local complete intersection of codimension $b$ inside $\P_m$.
\end{prop}

\begin{proof}
Denote with $\H_m=\PP_B\left(W\ot\EE\right)$ and $\H_m^{\circ}$ the open subset consisting of pairs $\left(z,\left[T\right]\right)$ such that $T^{\flat}:W\d\to\EE_z$ is surjective. The tautological inclusion $\OO_{\H_m}\left(-1\right)\su W\ot\EE$ induces a morphism of vector bundles on $\H_m$
\[
\lambda:W\d\ot\Q\d\longrightarrow\T\d\ot\OO_{\H_m}\left(1\right),
\]
whose fibre at $\left(z,\left[T\right]\right)$ is, up to a scalar, the composition $W\d\ot\Q\d\xrightarrow{T^{\flat}\ot\Q\d}\EE_z\ot\Q\d\to\T\d$. It is surjective on $\H_m\c$, so that $\K_m=\ker\lambda_{\vert\H_m\c}$ is a vector bundle of rank $\rk\K_m=m\left(N-k\right)-b>0$ and, in view of the definition of $\C_m$ as $Z\left(s_m\right)$,
\[
\C_{m\vert\H_m\c}\simeq\PP_{\H_m\c}\left(\K_m\right) ,
\]
which is nonempty; in particular $\C_{m\vert\H_m\c}$ is a smooth irreducible variety of dimension
\[
\dim\C_{m\vert\H_m\c}=\left(b+mk-1\right)+\left(m\left(N-k\right)-b\right)-1=mN-2.
\]
The inclusions $\EE\su S_r\ot\OO_B$ and $\Q\d\su S_r\d\ot\OO_B$ determine a morphism $\P_m\to\PP\left(\mathcal{M}_m\right)\times\PP\left(\mathcal{M}_m\d\right)$ that restricts to $\Phi$. Consider now $\Phi_{\vert\C_m\c}$. If $z\in B\c$ and $\rk T^{\flat}=k$ then $L=\EE_z\in\Sigma_{\mathrm{reg}}$, $\T_z=T_L \Sigma$ and $s_m  \left(z,\left[T\right],\left[Y\right]\right) = \kappa_L \left(T\ot Y\right)$, so that $\left(\left[T\right],\left[Y\right]\right)\in\N\c_{X_m}$ if and only if $s_m \left(z,\left[T\right],\left[Y\right]\right) = 0$ in view of  Lemma \ref{lem:CriterioConormale}. Moreover $L$ is reconstructed from $T$ as $L=\mathrm{Im}\left(T^{\flat}\right)$, so the inverse of $\Phi_{\vert\C_m\c}$ is $\left(\left[T\right],\left[Y\right]\right)\mapsto\left(\left(\nu_{\vert B\c}\right)^{-1}\left(\gamma\left(\left[T\right]\right)\right),\left[T\right],\left[Y\right]\right)$ with $\gamma$ as in Proposition \ref{prop:XcircLiscio}, hence $\Phi$ induces an isomorphism $\C_m\c\simeq\N\c_{X_m}$. Let now
\[
\D_m=\left\{ \left(z,\left[T\right]\right)\in\H_m\,\vert\,\rk T^{\flat}\le k-1\right\} ,
\]
this is fibrewise the determinantal variety of $k\times m$ matrices of rank at most $k-1$ hence, see for example \cite[Proposition 12.2]{Ha92}, it has codimension
\[
\mathrm{codim}_{\H_m}\D_m=m-k+1
\]
and the same holds for its preimage $\tilde{\D}_m\su\P_m$. On the other hand $\C_m$ is cut by the $b$ components of $s_m$ hence every irreducible component $Z$ of $\C_m$ satisfies $\mathrm{codim}_{\P_m}Z\le b$, but if $Z$ were contained in $\tilde{\D}_m$ then $\mathrm{codim}_{\P_m}Z\ge m-k+1$ which contradicts the assumption $m\ge k+b$. Hence every irreducible component $Z$ of $\C_m$ meets the open subset $\C_{m\vert\H_m\c}$, so that $Z\cap\C_{m\vert\H_m\c}$ is dense in $Z$ and $Z\su\overline{\C_{m\vert\H_m\c}}$; the latter being irreducible, $\C_m$ has a single irreducible component and  $\C_m = \overline{\C_{m\vert\H_m\c}}$. Moreover $\C_{m\vert\H_m\c}\to B$ is surjective since $m\ge k$, hence the preimage of the dense open subset $B\c\su B$, which is precisely $\C_m\c$, is a nonempty open subset of $\C_{m\vert\H_m\c}$ and therefore dense in it. We conclude that $\C_m=\overline{\C_m\c}$, which also proves $\left(2\right)$. In view of properness of $\Phi$ and Lemma \ref{lem:ConormaleSuXcirc}, it follows that 
\[
\Phi\left(\C_m\right)=\overline{\Phi\left(\C\c_m\right)}=\overline{\N_{X_m}\c}=\N_{X_m}
\]
therefore $\Phi:\C_m\to\N_{X_m}$ is proper and birational. Note that 
\[
\dim\P_m=b+\left(mk-1\right)+\left(m\left(N-k\right)-1\right)=mN+b-2=\dim\C_m+b
\]
so that $s_m$ is a regular section; since $\P_m$ is smooth, $\C_m$ is a local complete intersection, hence Cohen-Macaulay. Being generically smooth it satisfies $\left(R_{0}\right)$ and $\left(S_{1}\right)$, so it is reduced \cite[p. 71]{BH}, and being irreducible it is integral.
\end{proof}

\section{Stable polynomiality in $m$}\label{sec:polinomialita}

\subsection{Reduction to an integral over $B$}%~\\
Since $\C_m$ is the zero locus of a regular section of a bundle of rank $b$ over the smooth variety $\P_m$, we have
\begin{equation}\label{eq:ClasseCm}
\left[\C_m\right]=c_b\left(\T\d\ot\OO_{\P_m}\left(1,1\right)\right)\cap\left[\P_m\right]
\end{equation}
in $A^{\ast}\left(\P_m\right)$, see \cite[Example 3.2.16]{Ful}. Write $x=c_{1}\left(\OO_{\P_m}\left(1,0\right)\right)$ and $y=c_{1}\left(\OO_{\P_m}\left(0,1\right)\right)$ for the two relative hyperplane classes, and $\ga,\gb$ for the hyperplane classes of $\PP\left(\M_m\right)$ and $\PP\left(\M_m\right)\d$; by construction $\Phi \inv\ga=x$ and $\Phi \inv\gb=y$.

\begin{lem}\label{lem:IntegraleCmeNXm}
For $0\le j\le mN-2$ one has
\begin{equation}\label{eq:IntegraleCmeNXm}
\delta_j\left(X_m\right)=\int_{\N_{X_m}}\ga^{j}\gb^{mN-2-j}=\int_{\C_m}x^{j}y^{mN-2-j}.
\end{equation}
\end{lem}

\begin{proof}
In view of Proposition \ref{prop:ProprietaCm} the morphism $\Phi :\C_m\to\N_{X_m}$ is proper and birational between integral varieties, so $\Phi_{\ast}\left[\C_m\right]=\left[\N_{X_m}\right]$ and the second equality follows from the projection formula. For the first, pair \eqref{eq:DefinizioneGradiPolariConormale} with $\ga^{j}\gb^{mN-2-j}$: in the ambient $\PP\left(\M_m\right)\times\PP\left(\M_m\right)\d$ the term of index $i$ contributes only for $i=j$, the exponent of $\ga$ being at least $mN$ for $i<j$ and that of $\gb$ being at least $mN$ for $i>j$, and $\int\ga^{mN-1}\gb^{mN-1}=1$.
\end{proof}
 
Summing \eqref{eq:IntegraleCmeNXm} over $j$ and using \eqref{eq:ClasseCm} together with 
\[
c_b\left(\T\d\ot\OO\left(1,1\right)\right)=\sum_{i=0}^{b}c_i\left(\T\d\right)\left(x+y\right)^{b-i},
\]
see \cite[Example 3.2.2]{Ful}, we get
\begin{equation}\label{eq:gEDGrezza}
\ged\left(X_m\right)=\int_{\P_m}\Big[\sum_{i=0}^{b}c_i\left(\T\d\right)\left(x+y\right)^{b-i}\Big]
\cdot\Big[\sum_{j=0}^{mN-2}x^{j}y^{mN-2-j}\Big].
\end{equation}
Let $\rho_m:\P_m\to B$ be the structure map. Since $\rk\left(W\ot\EE\right)=mk$ and $\rk\left(W\d\ot\Q\d\right)=m\left(N-k\right)$, the projection formula for the two projective bundles gives
\begin{equation}\label{eq:FormulaPushforwardSegre}
\rho_{m\ast}\left(x^{mk-1+u}y^{m\left(N-k\right)-1+v}\right)
=s_u\left(W\ot\EE\right)s_v\left(W\d\ot\Q\d\right),\qquad u,v\ge0,
\end{equation}
all other monomials of the same total degree being sent to zero; here $s_{\bullet}$ denotes Segre classes, see \cite[\S3.1]{Ful}.

\begin{lem}\label{lem:Coefficiente2g}
For $g \ge0$ suppose
\begin{equation}\label{eq:CondizioneCombinatoria}
g\le mk-1\qquad\text{and}\qquad g\le m\left(N-k\right)-1 .
\end{equation}
Then
\[
\rho_{m\ast}\Big[\left(x+y\right)^{g}\sum_{j=0}^{mN-2}x^{j}y^{mN-2-j}\Big]
=2^{g}\sum_{\substack{u,v\ge0\\u+v=g}}s_u\left(W\ot\EE\right)s_v\left(W\d\ot\Q\d\right).
\]
\end{lem}

\begin{proof}
Expanding,
\[
\left(x+y\right)^{g}\sum_{j=0}^{mN-2}x^{j}y^{mN-2-j}
=\sum_{a=0}^{g}\binom ga\sum_{j=0}^{mN-2}x^{a+j}y^{mN-2-j+g-a},
\]
and the monomial $x^{mk-1+u}y^{m\left(N-k\right)-1+v}$ with $u+v=g$ occurs, for a given $a$, exactly when $j=mk-1+u-a$. It therefore occurs for every $a$ provided $0\le mk-1+u-a\le mN-2$ for all $0\le a\le g$: since $0\le u\le g$ and $0\le a\le g$ we have $-g\le u-a\le g$, hence $mk-1-g\le mk-1+u-a\le mk-1+g$, and \eqref{eq:CondizioneCombinatoria} gives $mk-1-g\ge0$ and $mk-1+g\le mN-2$. The coefficient is then $\sum_{a=0}^{g}\binom ga=2^{g}$, and \eqref{eq:FormulaPushforwardSegre} concludes.
\end{proof}

Applying Lemma \ref{lem:Coefficiente2g} with $g=b-i$ to each summand of \eqref{eq:gEDGrezza}, and using $W\ot\EE\simeq\EE^{\op m}$ and $W\d\ot\Q\d\simeq\left(\Q\d\right)^{\op m}$, we obtain the formula on which the rest of the paper rests:
\begin{equation}\label{eq:gEDFormulaBuona}
\ged\left(X_m\right)=\sum_{i=0}^{b}2^{b-i}\sum_{\substack{u,v\ge0\\u+v=b-i}}
\int_Bc_i\left(\T\d\right)\,s_u\left(\EE^{\op m}\right)\,s_v\left(\left(\Q\d\right)^{\op m}\right).
\end{equation}

\subsection{The polynomial $P_{n,k,r}$}%~\\
\begin{lem}\label{lem:GradoMClassiDiSegre}
Let $\F$ be a vector bundle. Then $m\mapsto s_u\left(\F^{\op m}\right)$ is a polynomial in $m$ of degree at most $u$.
\end{lem}

\begin{proof}
By Whitney's formula $s\left(\F^{\op m}\right) = c\left(\F\right)^{-m} = s\left(\F\right)^{m}$, so that the coefficient of $t^{u}$ is 
\[
\sum_{k_{0}+\dots+k_u=m}\binom{m}{k_{0},\dots,k_u}s_{1}\left(\F\right)^{k_{1}}\cdots s_u\left(\F\right)^{k_u},
\] 
with $k_{1}+2k_{2}+\dots+uk_u=u$. We conclude since the multinomial coefficient is a polynomial in $m$ of degree 
\[
k_{1}+\dots+k_u\le k_{1}+2k_{2}+\dots+uk_u=u .
\]
\end{proof}

\begin{thm}\label{thm:PolinomialitaStabileInM}
Fix $n,k,r$ with $1\le k<N=\binom{n+r-1}{r}$ and let $b=\dim\left(\Sigma\right)$. There is a polynomial $P=P_{n,k,r}\in\QQ\left[t\right]$ such that
\[
\ged\left(X_m\right)=P\left(m\right)\qquad\text{for every }m\ge k+b .
\]
Its degree is exactly $b$ and its leading coefficient is $\frac{4^{b}}{b!}\deg\left(\Sigma_{k,r}^{n}\right)$, where $\deg$ denotes the degree in the Plücker embedding of $\Gr_{k}\left(S_{r}\right)$; in particular
\[
\ged\left(X_m\right)=\frac{4^{b}}{b!}\,\deg\left(\Sigma_{k,r}^{n}\right)m^{b}+O\left(m^{b-1}\right) .
\]
\end{thm}

\begin{proof}
For $m\ge k+b$ Proposition \ref{prop:ProprietaCm} applies, and \eqref{eq:CondizioneCombinatoria} holds for every $g\le b$ because 
\[
mk\ge\left(k+b\right)k\ge k+b>b \quad\text{and} \quad m\left(N-k\right)\ge m\ge k+b>b ,
\]
hence \eqref{eq:gEDFormulaBuona} holds. In it $B,\EE,\Q,\T$ and $b$ do not depend on $m$, while in view of Lemma \ref{lem:GradoMClassiDiSegre} each product $s_u\left(\EE^{\op m}\right)s_v\left(\left(\Q\d\right)^{\op m}\right)$ is polynomial in $m$ of degree at most $u+v=b-i\le b$; this shows that $\ged\left(X_m\right)$ itself is a polynomial in $m$ of degree at most $b$. Let us now consider the leading coefficient of $m^{b}$. A term of $m$-degree $b$ requires $u+v=b$, hence $i=0$. Expanding $s\left(\EE^{\op m}\right)=c\left(\EE\right)^{-m}$ and $s\left(\left(\Q\d\right)^{\op m}\right)=c\left(\Q\d\right)^{-m}$ as exponentials, the coefficient of $m^{u}$ in $s_{u}\left(\EE^{\op m}\right)$ is $\frac{1}{u!}\left(-c_{1}\left(\EE\right)\right)^{u}$ and that of $m^{v}$ in $s_{v}\left(\left(\Q\d\right)^{\op m}\right)$ is $\frac{1}{v!}\left(-c_{1}\left(\Q\d\right)\right)^{v}=\frac{1}{v!}c_{1}\left(\Q\right)^{v}$. Since $c_{1}\left(\EE\right)+c_{1}\left(\Q\right)=0$, the coefficient of $m^{b}$ is 
\[
2^{b}\sum_{u+v=b}\frac{1}{u!\,v!}\int_{B}c_{1}\left(\Q\right)^{b} = \frac{4^{b}}{b!}\int_{B}c_{1}\left(\Q\right)^{b} .
\]
Now $c_{1}\left(\Q\right)$ denotes the pullback of the Plücker hyperplane class, $\nu$ is birational and $\dim\left(\Sigma_{k,r}^{n}\right)=b$, so the integral is $\deg \left( \Sigma_{k,r}^{n} \right)>0$, in particular $\deg P=b$.
\end{proof}

\begin{rem}\label{rem:DueCondizioni}
The hypothesis $m\ge k+b$ enters the proof three times and for unrelated reasons: through Lemma \ref{lem:ConormaleSuXcirc}, to ensure $m\left(N-k\right)>b$, that is $\rk\K_{m}>0$; through Proposition \ref{prop:ProprietaCm}, to ensure that no irreducible component of $\C_{m}$ is contained in the locus where $\rk T^{\flat}<k$; and through \eqref{eq:CondizioneCombinatoria}, to ensure that the coefficient of Lemma \ref{lem:Coefficiente2g} is $2^{g}$ rather than a truncated binomial sum. Each of the three may hold under weaker assumptions on $m$, and they need not be binding at the same time: for $k=1$ the second is vacuous, every nonzero $T\in W\ot\EE_{z}$ having $\rk T^{\flat}=1$, and the first is implied by the third, so that the bound $m\ge n=k+b$ comes from the combinatorics alone; for $n=k=r=2$ all three are weaker than $m\ge k+b=4$. The examples of Appendix \ref{sec:esempi} show that the bound is nevertheless sharp as a uniform statement.
\end{rem}

\section{Joint polynomiality in $m$ and $r$}
\label{sec:polinomialitaCongiunta}

Theorem \ref{thm:PolinomialitaStabileInM} holds for a fixed activation exponent $r$. We now let $r$ vary as well. The point, as we are going to see, is that for $r\ge k$ the base of the incidence construction stabilises, so that a single $B$ serves for all $r$ and the whole dependence on $r$ is carried by one bundle, whose Chern classes turn out to be polynomial in $r$.

\subsection{Stabilisation of the base}\label{sec:baseFissa}%~\\
We  denote by $\Hilb_{k}\PP\left(V\d\right)$ the Hilbert scheme of length $k$ subschemes of  $\PP\left(V\d\right)$, by $p:\mathscr{Z} \to\Hilb_{k} \PP\left(V\d\right)$ its universal family and by $\xi: \mathscr{Z} \to\PP\left(V\d\right)$ the projection, so that $p$ is  finite locally free of degree $k$; see \cite[Theorem 1.1 and Corollary 1.2]{HS04} and \cite[Tag 0B94]{Stacks}.

Since $p$ is finite locally free of degree $k$, the sheaf $\EE_{r}:=p_{\ast}\xi\inv\OO\left(r\right)$ is locally free of rank $k$, and in view of Lemma \ref{lem:Regolarita} the evaluation $S_{r}\d\ot\OO\to\EE_{r}$ is surjective for $r\ge k-1$. Dualising gives a subbundle $\EE_{r}\d\su S_{r}\ot\OO$, that is a morphism
\[
\eta_{r}:\Hilb_{k}\PP\left(V\d\right)\to\Gr_{k}\left(S_{r}\right),\qquad
\left[Z\right]\mapsto\left(I_{Z}\right)_{r}^{\perp},
\]
along which $\EE_{r}\d$ is the pullback of the tautological subbundle.

\begin{prop}\label{prop:HilbertImmersione}
For $r\ge k$ the morphism $\eta_{r}$ is a closed embedding.
\end{prop}

\begin{proof}
The statement for $r=k$ is classical: by the proof of Lemma \ref{lem:Regolarita}, the Gotzmann number of the constant Hilbert polynomial $k$ is $k$, and $\left\{ k\right\}$ is a supportive set of degrees, so that \cite[Lemma 4.1, Proposition 4.2 and (21)]{HS04} apply. 
%The references we are aware of, such as \cite{HS04} and \cite{BLR13}, realise the Hilbert scheme inside a Grassmannian in that single degree only; we therefore include the passage from $r=k$ to $r>k$.
Let now $r>k$, let $\mathfrak{m}=\Sym^{\ge1}V$ be the irrelevant ideal, let $\mathscr{Z}\su\PP\left(V\d\right)\times T$ be a family of length $k$ subschemes over an arbitrary base and let $p:\mathscr{Z}\to T$ be the projection. Since $p$ is finite and flat of degree $k$, the sheaves $\G_{d}=p_{\ast}\OO_{\mathscr{Z}}\left(d\right)$ are locally free of rank $k$ and their formation commutes with base change; in view of Lemma \ref{lem:Regolarita} the evaluation $S_{d}\d\ot\OO_{T}\to\G_{d}$ is surjective for $d\ge k-1$, surjectivity being checked on the fibres, so that $\left(I_{\mathscr{Z}}\right)_{d}$ is a subbundle of $S_{d}\d\ot\OO_{T}$ with quotient $\G_{d}$.

Since $S_{r-k}\d\cdot\left(I_{\mathscr{Z}}\right)_{k}\su\left(I_{\mathscr{Z}}\right)_{r}$, we have a morphism
\[
\G_{k}\longrightarrow\Hom\left(S_{r-k}\d,\G_{r}\right),\qquad
\overline{f}\longmapsto\left(g\mapsto\overline{gf}\right),
\]
which is injective on every fibre: if $Z$ has length $k$ and $f\in S_{k}\d$ satisfies $S_{r-k}\d\cdot f\su\left(I_{Z}\right)_{r}$, then $\mathfrak{m}^{r-k}f\su I_{Z}$, whence $f\in\left(I_{Z}:\mathfrak{m}^{r-k}\right)=I_{Z}$, the ideal of a subscheme being  saturated. It follows that it is an inclusion of subbundles, hence
\[
\left(I_{\mathscr{Z}}\right)_{k}
=\ker\Big(S_{k}\d\ot\OO_{T}\longrightarrow\Hom\left(S_{r-k}\d,\G_{r}\right)\Big),
\]
is a subbundle depending only on $\left(I_{\mathscr{Z}}\right)_{r}$. It follows that $\eta_{k}$ factors through $\eta_{r}$ on $T$-points for every $T$; since $\eta_{k}$ is a monomorphism, so is $\eta_{r}$ and we conclude since a proper monomorphism is a closed embedding.
\end{proof}

The locus $\Hilb_{k}\c\PP\left(V\d\right)$ of reduced subschemes is open, being the étale  locus of $p$, and is smooth irreducible of dimension $k\left(n-1\right)$, being the quotient of the complement of the diagonals in $\PP\left(V\d\right)^{k}$ by the symmetric group; we write $\Hilb_{k}^{\mathrm{sm}}\PP\left(V\d\right)$ for the reduced closure, called the smoothable component, which is irreducible of the same dimension and in general singular. Note that every subscheme of length $2$ is smoothable, so that
\[
\Hilb_{2}^{\mathrm{sm}}\PP\left(V\d\right)=\Hilb_{2}\PP\left(V\d\right),
\]
which is smooth and irreducible: for $n\ge3$ see \cite[\S2.1]{BFR19}, while $\Hilb_{2}\PP^{1}\simeq\PP^{2}$.

The following result is the one that states the independence of $\Sigma_{k,r}^{n}$ on $r$.

\begin{prop}\label{prop:BaseFissa}
Let $r\ge k$. Then $\eta_{r}$ restricts to an isomorphism
\[
\Hilb_{k}^{\mathrm{sm}}\PP\left(V\d\right)\overset{\sim}{\longrightarrow}\Sigma_{k,r}^{n}
\]
under which $\EE_{r}\d$ corresponds to the tautological bundle $\UU$ of Section \ref{sec:neurovarieta}. 
\end{prop}

\begin{proof}
Let $Z=\left\{ \left[\ell_{1}\right],\dots,\left[\ell_{k}\right]\right\} $ be reduced. The
pairing of Section \ref{sec:neurovarieta} gives $\left\langle f,\ell_{i}^{r}\right\rangle =f\left(\ell_{i}\right)$, so that $\left\langle \ell_{1}^{r},\dots,\ell_{k}^{r}\right\rangle
\su\left(I_{Z}\right)_{r}^{\perp}$; by Lemma \ref{lem:Regolarita} both sides have dimension
$k$, whence equality. Therefore $\eta_{r}\left(\Hilb_{k}\c\PP\left(V\d\right)\right)$ is the locus whose closure is $\Sigma_{k,r}^{n}$ by definition. By Proposition \ref{prop:HilbertImmersione} the image of $\Hilb_{k}^{\mathrm{sm}}\PP\left(V\d\right)$ is closed and irreducible and contains that locus as a dense subset, hence equals $\Sigma_{k,r}^{n}$; both schemes being reduced, $\eta_{r}$ is an isomorphism onto it, and it carries $\EE_{r}\d$ to $\UU$ because $\EE_{r}\d$ is the pullback of the tautological subbundle. The last assertion follows since neither $\Hilb_{k}^{\mathrm{sm}}\PP\left(V\d\right)$ nor its universal family involve $r$.
\end{proof}

It follows from Proposition \ref{prop:BaseFissa} that there are a variety $\Sigma$, a finite locally free family $p:\mathscr{Z}\to\Sigma$ of degree $k$ and a morphism $\xi: \mathscr{Z} \to \PP \left(V\d\right)$, none of which depends on $r$, such that
\[
\UU\simeq\EE_{r}\d,\qquad\EE_{r}=p_{\ast}\xi\inv\OO\left(r\right),
\]
for every $r\ge k$: the whole dependence on $r$ of the construction of Section \ref{sec:conormale} is concentrated in this single bundle. 
%In particular Proposition \ref{prop:BaseFissa} reproves Proposition \ref{prop:DimensioneSigma} for $r\ge k$. 
The Nash blow-up of $\Sigma$ and its tautological bundle are intrinsic to $\Sigma$, so that the pair $\left(B,\T\right)$ of Section \ref{sec:conormale} is of the type considered in \cite[Example 4.2.9]{Ful}: a proper birational $B\to\Sigma$ along which $\T\d$ is a locally free rank $b$ quotient of $\nu\inv\Omega_{\Sigma}$. Since $\EE$ and $\Q$ are pulled back from $\Sigma$, the projection formula shows that the integrals in \eqref{eq:gEDFormulaBuona} depend on $B$ only through the classes $\nu_{\ast}\left(c_{i}\left(\T\d\right)\cap\left[B\right]\right)$, which up to sign are the Chern-Mather classes of $\Sigma$ and are therefore independent of the choice of $B$ \cite[Example 4.2.9(b)]{Ful}. As $\Sigma$ does not depend on $r$, one and the same $B$ and $\T$ serve for every $r\ge k$. By contrast the embedding $\T\su\EE\d\ot\Q$ used to define $s_{m}$ does depend on $r$, since $\rk\Q=N-k$ grows with $r$; only the classes $c_{i}\left(\T\d\right)$ do not.

\begin{rem}\label{rem:SogliaRK}
The hypothesis $r\ge k$ is necessary. For $r=1$ one has $\ell^{1}=\ell$, so that $\Sigma_{k,1}^{n}=\Gr_{k}\left(V\d\right)$, of dimension $k\left(n-k\right)$, whereas $\Hilb_{k}^{\mathrm{sm}} \PP\left(V\d\right)$ has dimension $k\left(n-1\right)$ for every $r$. For $2\le k<n$ the two differ, hence Proposition \ref{prop:BaseFissa} fails; and since
$b=\dim\left(\Sigma\right) = \deg P_{n,k,r}$ for Theorem \ref{thm:PolinomialitaStabileInM}, no base uniform in $r$ can exist.
\end{rem}

\subsection{The polynomial $Q_{n,k}$}\label{sec:teoremaCongiunto}%~\\
In view of Proposition \ref{prop:BaseFissa} the only object entering \eqref{eq:gEDFormulaBuona} which still depends on $r$ is the bundle $\EE_{r} = p_{\ast} \xi \inv\OO\left(r\right)$. We show now that its Chern character is polynomial in $r$. Indeed, in view of a more general version of the Riemann-Roch Theorem, in the form valid for singular schemes \cite[Chapter 18]{Ful}, the twist by $\OO\left(r\right)$ turns into multiplication by $e^{rH}$, and capping the latter with a class living in dimensions at most $b$ truncates the exponential to a polynomial, the degree in $r$ of each graded piece being controlled by its codimension. Theorem \ref{thm:PolinomialitaCongiunta} then follows by combining this with the polynomiality in $m$ of Theorem \ref{thm:PolinomialitaStabileInM}.

\begin{lem}\label{lem:ChernPolinomialiInR}
Let $r\ge k$. For every $i$ the class $c_{i}\left(\EE\right)\in A^{i}\left(B\right)_{\QQ}$ is given by a polynomial in $r$ of degree at most $i$, whose coefficients do not depend on $r$. The same holds for $s_{u}\left(\EE\right)$, $c_{v}\left(\Q\d\right)$ and $s_{v}\left(\Q\d\right)$, of degree at most $u$, respectively $v$.
\end{lem}

\begin{proof}
By Proposition \ref{prop:BaseFissa} we may pull back the universal family to $B$ and regard $p:\mathscr{Z}\to B$ as finite locally free of degree $k$, with $\xi:\mathscr{Z}\to \PP\left( V\d\right)$ and $\EE\simeq\EE_{r}\d$, $\EE_{r}=p_{\ast}\xi\inv\OO\left(r\right)$; none of $\mathscr{Z}$, $p$, $\xi$ depends on $r$. Put $H=c_{1}\left(\xi\inv\OO\left(1\right)\right)$. Since $p$ is finite locally free, $\EE_{r}$ is locally free and $p_{\ast}\left[\xi\inv\OO\left(r\right)\right]=\left[\EE_{r}\right]$ in $K_{0}\left(B\right)$, the higher direct images vanishing, so \cite[Corollary 18.3.1(c)]{Ful} applies and gives
\[
p_{\ast}\left(e^{rH}\cap\Td\left(\mathscr{Z}\right)\right)
=\ch\left(\EE_{r}\right)\cap\Td\left(B\right)
=\ch\left(\EE_{r}\right)\cap\td\left(T_{B}\right)\cap\left[B\right],
\]
the second equality by \cite[Corollary 18.3.1(b)]{Ful}, $B$ being smooth. Note that $\mathscr{Z}$ is singular in general, so that $\td\left(T_{\mathscr{Z}}\right)$ is not available and the singular form of Riemann-Roch is needed here, and $\Td\left(\mathscr{Z}\right)$ denotes $\tau_{\mathscr{Z}}\left(\OO_{\mathscr{Z}}\right)$, the Todd class of an arbitrary scheme in the sense of \cite[Chapter 18]{Ful}, not the usual $\td\left(T_{\mathscr{Z}}\right)$. We point out that the class $\tau:=\Td\left(\mathscr{Z}\right)\in A_{\ast}\left(\mathscr{Z}\right)_{\QQ}$ does not depend on $r$.

As $p$ is finite, $\dim\mathscr{Z}=\dim B=b$, hence $\tau\in\bigoplus_{j\le b}A_{j}\left(\mathscr{Z}\right)_{\QQ}$; capping with $H$ lowers dimension by one, so $H^{a}\cap\tau=0$ for $a>b$ and
\[
e^{rH}\cap\tau=\sum_{a=0}^{b}\frac{r^{a}}{a!}\,H^{a}\cap\tau .
\]
Since $p$ preserves dimensions, the contribution of the $a$-th summand to  $A_{b-c}\left(B\right)_{\QQ}$ is $p_{\ast}\left(H^{a}\cap\tau_{b-c+a}\right)$, which vanishes unless $a\le c$. Comparing the components of dimension $b-c$ on the two sides and arguing by induction on $c$, using that $\td\left(T_{B}\right)\cap\left[B\right]$ has leading term $\left[B\right]$, we conclude that $\ch_{g}\left(\EE_{r}\right)$ is a polynomial in $r$ of degree at most $g$.

Finally $c_{i}$ is a universal polynomial with rational coefficients in $\ch_{1},\dots,\ch_{i}$, isobaric of weight $i$, so $\deg_{r}c_{i}\left(\EE_{r}\right)\le i$; one has $c_{i}\left(\EE\right)=\left(-1\right)^{i}c_{i}\left(\EE_{r}\right)$, the Segre classes $s_{u}\left(\EE\right)$ are universal polynomials in $c_{1}\left(\EE\right),\dots,c_{u}\left(\EE\right)$ isobaric of weight $u$, and $c\left(\Q\right)=c\left(\EE\right)^{-1}$ by Whitney's formula and $c_i\left(\Q\d\right)=(-1)^i c_i\left(\Q\right)$; since Segre classes are universal isobaric polynomials in the Chern classes, all the remaining assertions follow.
\end{proof}

\begin{thm}\label{thm:PolinomialitaCongiunta}
There is a polynomial $Q=Q_{n,k}\in\QQ\left[m,r\right]$, with $\deg_{m}Q\le b$ and $\deg_{r}Q\le b$, such that
\[
\ged\left(X_{k,r}^{n,m}\right)=Q\left(m,r\right)\qquad\text{for all }r\ge k,\ m\ge kn .
\]
\end{thm}

\begin{proof}
In view of Proposition \ref{prop:DimensioneSigma}, for $r\ge k$ one has $b=k\left(n-1\right)$, so the threshold $m\ge k+b=kn$ of Theorem \ref{thm:PolinomialitaStabileInM} does not depend on $r$, moreover $k<N=\binom{n+r-1}r$ holds since $N\ge r+1>k$. Hence \eqref{eq:gEDFormulaBuona} is available on the whole quadrant. In each summand
\[
T_{i,u,v}=2^{\,b-i}\int_{B}c_{i}\left(\T\d\right)\,s_{u}\left(\EE^{\op m}\right)\,
s_{v}\left(\left(\Q\d\right)^{\op m}\right),\qquad u+v=b-i,
\]
the class $c_{i}\left(\T\d\right)$ is independent of both $m$ and $r$ by Proposition \ref{prop:BaseFissa}. In view of Lemma \ref{lem:GradoMClassiDiSegre}, $s_{u}\left(\EE^{\op m}\right)$ is a polynomial in $m$ of degree at most $u$ whose coefficients are the monomials $\prod_{j}s_{j}\left(\EE\right)^{k_{j}}$ with $\sum_{j}jk_{j}=u$; each such monomial has $r$-degree at most $\sum_{j}jk_{j}=u$ by Lemma \ref{lem:ChernPolinomialiInR}. Hence $s_{u}\left(\EE^{\op m}\right)$ is polynomial in $\left(m,r\right)$ with $\deg_{m}\le u$ and $\deg_{r}\le u$, and likewise $s_{v}\left(\left(\Q\d\right)^{\op m}\right)$ with $\deg_{m}\le v$ and $\deg_{r}\le v$. It follows that $T_{i,u,v}$ is a polynomial in $\left(m,r\right)$ with $\deg_{m},\deg_{r}\le u+v=b-i\le b$. Summing over $i,u,v$ we conclude.
\end{proof}

We evaluate formula \eqref{eq:gEDFormulaBuona} in closed form in two cases in Appendix \ref{sec:esempi}: for $k=1$, where it gives the polynomial $P_{n,r}$ of \eqref{eq:PolinomioCasoK1} whose specialisation at $r=1$ reproduces the values tabulated in \cite[Table 2]{KMQS25}; and for $n=k=2$, where Proposition \ref{prop:FormulaCasoBinario} gives
\[
\ged\left(X_{2,r}^{2,m}\right)=8\left(r-1\right)^{2}m^{2}-12\left(r-1\right)m+3 ,
\qquad m\ge4 .
\]
The slice $r=2$ is the determinantal case of \cite[Theorem 6.5]{KLW24}, whose stated range of validity we sharpen in Remark \ref{rem:ConfrontoKLW}.

\section{A general statement and a glimpse of deep networks}\label{sec:Finale}

Note that nothing in Sections \ref{sec:conormale} and \ref{sec:polinomialita} uses the explicit description of $\Sigma$ as a variety of spaces of powers. The only input from Section \ref{sec:neurovarieta} is the presentation of $X_{m}$ as the union of the projective subspaces $\PP\left(W\ot L\right)$, $L\in\Sigma$, and this makes the whole proof of the polynomiality in $m$ available for an arbitrary $\Sigma$.

\begin{thm}\label{thm:GeneraleSigma}
Let $S$ be a complex vector space of dimension $N$, let $1\le k<N$ and let $\Sigma\su\Gr_{k}\left(S\right)$ be an irreducible closed subvariety of dimension $b$. For a vector space $W$ of dimension $m$ set
\[
X_{\Sigma,m}=\bigcup_{L\in\Sigma}\PP\left(W\ot L\right)\su\PP\left(W\ot S\right).
\]
Then $X_{\Sigma,m}$ is an irreducible projective variety, of dimension $km+b-1$ for every $m\ge k$, and there is a polynomial $P_{\Sigma}\in\QQ\left[t\right]$ such that
\[
\ged\left(X_{\Sigma,m}\right)=P_{\Sigma}\left(m\right)\qquad\text{for every }m\ge k+b .
\]
Its degree is exactly $b$ and its leading coefficient is $\frac{4^{b}}{b!}\deg\Sigma$, the degree being taken in the Plücker embedding of $\Gr_{k}\left(S\right)$.
\end{thm}

\begin{proof}
Let $\UU$ be the restriction to $\Sigma$ of the tautological subbundle, let $\Y=\PP_{\Sigma}\left(W\ot\UU\right)$ and let $\mu:\Y\to\PP\left(W\ot S\right)$ be the second projection; then $X_{\Sigma,m}=\mu\left(\Y\right)$ is closed and irreducible, $\mu$ being proper and $\Y$ a projective bundle over an irreducible base. This is the conclusion of Lemma \ref{lem:IncidenzaPropria}, which here holds by definition; its second assertion holds as well, since $T\in W\ot L$ with $\rk T^{\flat}=k$ forces $\mathrm{Im}\left(T^{\flat}\right)=L$.

From this point on the arguments apply verbatim. The proofs of Proposition \ref{prop:XcircLiscio} and of Corollaries \ref{cor:ConoAffine} and \ref{cor:SpazioTangente} use only $\Y$, $\mu$ and the tangent spaces of $\Sigma$; those of Section \ref{sec:conormale} use only these together with the Nash blow-up of $\Sigma$ and a resolution $B$ of it, and the hypothesis $m\left(N-k\right)>b$ of Lemma \ref{lem:ConormaleSuXcirc} follows from $k<N$ and $m\ge k+b$ as in Theorem \ref{thm:PolinomialitaStabileInM}; those of Section \ref{sec:polinomialita} take place entirely on $B$.
\end{proof}

In this sense Theorem \ref{thm:PolinomialitaStabileInM} is the case $\Sigma=\Sigma_{k,r}^{n}$, while Theorem \ref{thm:PolinomialitaCongiunta} is specific to $\Sigma_{k,r}^{n}$, since such is the identification of Proposition \ref{prop:BaseFissa}.

The interest of this generality is that the same presentation is available for fully connected networks of arbitrary depth. Consider an architecture $\left(d_{0},\dots,d_{h}\right)$ with activation degree $r$ at every hidden layer, that is the family of maps
\[
f=\gt_{h}\circ\rho_{r}\circ\gt_{h-1}\circ\dots\circ\rho_{r}\circ\gt_{1}:V\to W ,
\qquad V=\RR^{d_{0}},\quad W=\RR^{d_{h}},
\]
homogeneous of degree $r^{h-1}$; put $k=d_{h-1}$ and $S=\Sym^{r^{h-1}}\left(V\d\right)$. Writing $g_{1},\dots,g_{k}\in S$ for the polynomials computed by the neurons of the last hidden layer and $w_{1},\dots,w_{k}\in W$ for the columns of $\gt_{h}$, one has $f=\sum_{j}w_{j}\ot g_{j}$, exactly as in Section \ref{sec:neurovarieta}. If the $g_{j}$ are linearly independent for a general choice of $\gt_{1},\dots,\gt_{h-1}$, let $\Sigma\p\su\Gr_{k}\left(S\right)$ be the closure of the locus of the spaces $\left\langle g_{1},\dots,g_{k}\right\rangle$ so obtained: it is irreducible since the parameters form an irreducible space. The argument of Lemma \ref{lem:IncidenzaPropria} then identifies the neurovariety of the architecture with $X_{\Sigma\p,d_{h}}$ and Theorem \ref{thm:GeneraleSigma} applies: for a fixed $\left(d_{0},\dots,d_{h-1}\right)$ and $r$, the generic $\ED$-degree of a deep polynomial network is a polynomial in the output dimension $d_{h}$, of degree $b\p=\dim\Sigma\p$ and leading coefficient $4^{b\p}\deg\left(\Sigma\p\right)/b\p!$, as soon as $d_{h}\ge d_{h-1}+b\p$.

What is missing is the geometry of $\Sigma\p$. Its dimension is accessible from the literature on the dimension of neurovarieties, since $b\p=\dim X_{\Sigma\p,d_{h}}-d_{h-1}d_{h}+1$ for $d_{h}\ge d_{h-1}$, compare \cite{KTB19,KLW24}; its degree in the Plücker embedding, which governs the leading coefficient, seems not to be known already for $h=3$. For $h=2$ one has $\Sigma\p=\Sigma_{k,r}^{n}$, whose dimension is Proposition \ref{prop:DimensioneSigma} and whose degree we compute for $k=2$ in Lemma \ref{lem:GradoSigma}.

\appendix

\section{A genericity statement for the isotropic quadric} \label{app:Genericity}

Over $\CC$, the statement proved below is the genericity assumption underlying Proposition \ref{prop:ProprietaGed}$\left(1\right)$: it is the hypothesis of \cite[Theorem 5.4]{DHOST}, and it is asserted for a general $q$, without proof, in the remark following \cite[Theorem 2.10]{KMQS25}. We include a proof because it is a plain dimension count on an incidence variety, which uses neither the ground field nor the characteristic beyond $\mathrm{char}\kk\neq2$; the restriction on the characteristic only serves to identify symmetric bilinear forms with quadratic forms and to make the associated polarity well behaved. 
%We stress that it is the geometric statement that generalizes: the $\ED$-degree itself is defined only over $\RR$ and $\CC$, since it requires a metric.

\begin{prop}\label{prop:GenericitaGenerale}
Let $\kk$ be an algebraically closed field with $\mathrm{char}\kk\neq2$, let $V=\kk^{M+1}$ and let $X\subsetneq\PP^M$ be a reduced and irreducible subvariety. Denote by $\mathcal{S}=\PP\left(\Sym^{2}V\d\right)$ the space of symmetric bilinear forms on $V$ up to scaling and, for $q\in\mathcal{S}$ nondegenerate, by
\[
\Delta_q=\overline{\left\{ \left(\left[x\right],\left[qx\right]\right)\,\vert\,x\in V\setminus\left\{ 0\right\} \right\} }\su\PP^M\times\PP\left(V\d\right)
\]
the diagonal induced by $q$. Then there exists a nonempty Zariski open subset $U\su\mathcal{S}$, consisting of nondegenerate forms, such that
\[
\N_X \cap\Delta_q=\e\qquad\text{for every }q\in U .
\]
\end{prop}

\begin{proof}
Consider the incidence variety
\[
\mathcal{Z}=\left\{ \left(\xi,\left[q\right]\right)\in\N_X\times\mathcal{S}\,\vert\,qx\in\left\langle y\right\rangle \right\} ,
\qquad \xi=\left(\left[x\right],\left[y\right]\right)\in\N_X \su\PP^M\times\PP\left(V\d\right),
\]
which is closed
%the condition $qx\in\left\langle y\right\rangle $ being bilinear in $q$ and $\left(x,y\right)$, 
and expresses that the polarity induced by $q$ maps $\left[x\right]$ to $\left[y\right]$ whenever $qx\neq0$. It comes with the two projections
\[
\xymatrix{ & \mathcal{Z}\ar[dl]_{\pi_{1}}\ar[dr]^{\pi_{2}}\\
\N_X & & \mathcal{S}
}
\]
We first compute $\dim\mathcal{Z}$ by means of $\pi_{1}$. Fix $\xi=\left(\left[x\right],\left[y\right]\right)\in\N_X $, so that $x\in V$ and $y\in V\d$ are nonzero, and consider the evaluation map
\[
\phi_x:\Sym^{2}V\d\to V\d,\qquad\phi_x\left(q\right)=qx .
\]
It is surjective: given $\ell\in V\d$, choose a basis of $V$ with $x=e_{0}$ and let $q$ be the symmetric bilinear form whose matrix has the components of $\ell$ as first row and first column and zero elsewhere; then $qx=\ell$. Composing with the quotient $V\d\to V\d/\left\langle y\right\rangle $, which is surjective onto an $M$-dimensional space, we obtain a surjection $\psi_{\xi}:\Sym^{2}V\d\to V\d/\left\langle y\right\rangle $ whose kernel has codimension $M$, and by construction $\pi_{1}^{-1}\left(\xi\right)=\PP\left(\ker\psi_{\xi}\right)$. Thus $\pi_{1}$ is surjective with fibres of constant dimension $\dim\mathcal{S}-M$, whence
\[
\dim\mathcal{Z}=\dim\N_X +\left(\dim\mathcal{S}-M\right) = \left(M-1\right)+\left(\dim\mathcal{S}-M\right)=\dim\mathcal{S}-1 ,
\]
where $\dim\N_X =M-1$ because $X$ is a proper subvariety of $\PP^M$.

Now consider $\pi_{2}$. Since $\mathcal{Z}$ is closed in $\N_X \times\mathcal{S}$ and $\N_X $ is projective, the image $\pi_{2}\left(\mathcal{Z}\right)$ is closed, and it is a proper closed subset of $\mathcal{S}$ because  $\dim\pi_{2}\left(\mathcal{Z}\right)\le\dim\mathcal{Z}<\dim\mathcal{S}$. Denoting by $D\su\mathcal{S}$ the hypersurface of degenerate forms, the set
\[
U=\mathcal{S}\setminus\left(\pi_{2}\left(\mathcal{Z}\right)\cup D\right)
\]
is open and nonempty, $\mathcal{S}$ being irreducible, and every $q\in U$ is nondegenerate and satisfies $\N_X \cap\Delta_q=\e$.
\end{proof}

\begin{rem}
For $\kk=\CC$ and $q\in U$, \cite[Theorem 5.4]{DHOST} gives $\ED_q\left(X\right)=\sum_i\delta_i\left(X\right)=\ged\left(X\right)$, which is Proposition \ref{prop:ProprietaGed}$\left(1\right)$.
\end{rem}

\section{Examples}\label{sec:esempi}

There are two cases in which \eqref{eq:gEDFormulaBuona} can be applied explicitly: $k=1$ where $\Sigma$ is a Veronese variety and $k=2$ with $n=2$ where $\Sigma\simeq\Hilb _2 \PP^1\simeq\PP^2$. In both, \eqref{eq:gEDFormulaBuona} can be evaluated in closed form and the answer compared with the literature; the comparison also shows how far the threshold $m\ge k+b$ of Theorem \ref{thm:PolinomialitaStabileInM} is from being optimal. For $k=2$ and arbitrary $n$ the base is still explicit, and smooth: it is $\Hilb_{2}\PP^{n-1}$, so that one may take $B=\Sigma$ and $\T=T_{\Sigma}$. We do not evaluate \eqref{eq:gEDFormulaBuona} in closed form there, but we compute the Plücker degree of $\Sigma_{2,r}^{n}$, which by Theorem \ref{thm:PolinomialitaStabileInM} determines the
leading term of $P_{n,2,r}$.

\subsection{The case $k=1$}\label{sec:esempiK1}%~\\
Here $X_{m}=\Seg\left(\PP\left(W_{m}\right)\times v_{r}\PP\left(V\d\right)\right)$ is smooth, so that $\ged\left(X_{m}\right)$ is computable by \eqref{eq:FormulaChernLiscio}, the only exception being $r=m=1$, where $X_{1}=\PP\left(\M_{1}\right)$ and $\ged$ is not defined. A closed formula is given in \cite[Theorem 3.5]{KMQS25}, $X_{m}$ being the Segre-Veronese variety $\mathcal{V}_{\left(1,r\right),\left(m-1,n-1\right)}$. The comparison with \eqref{eq:gEDFormulaBuona} therefore checks the whole construction against known values, and shows that the threshold cannot be lowered uniformly. Throughout, $\Sigma=v_{r}\PP\left(V\d\right)\simeq\PP^{n-1}$ is smooth of dimension $b=n-1$, one may take $B=B\c=\Sigma$, and $h$ denotes the hyperplane class of $\PP\left(V\d\right)$, so that $\int_{B}h^{n-1}=1$.

\begin{lem}\label{lem:FibratiCasoK1}
With the above identifications $\EE\simeq\OO\left(-r\right)$, $c\left(\Q\d\right)=\left(1+rh\right)^{-1}$ and $\T=T_{\PP^{n-1}}$; consequently
\[
s_{u}\left(\EE^{\op m}\right)=\binom{m+u-1}{u}r^{u}h^{u},\qquad
s_{v}\left(\left(\Q\d\right)^{\op m}\right)=\binom mv r^{v}h^{v},\qquad
c_{i}\left(\T\d\right)=\left(-1\right)^{i}\binom ni h^{i}.
\]
\end{lem}

\begin{proof}
The tautological subbundle of $\Gr_{1}\left(S_{r}\right)=\PP\left(S_{r}\right)$ is $\OO\left(-1\right)$ and $v_{r}\inv\OO\left(1\right)=\OO\left(r\right)$, whence $\EE\simeq\OO\left(-r\right)$ and $c\left(\EE\right)=1-rh$. From $0\to\EE\to S_{r}\ot\OO_{B}\to\Q\to0$ we get $c\left(\Q\right)=\left(1-rh\right)^{-1}$ and $c\left(\Q\d\right)=\left(1+rh\right)^{-1}$; and $\T=T_{\Sigma}$ because $\Sigma$ is smooth, so $c\left(\T\d\right)=\left(1-h\right)^{n}$. The three formulas follow from $s\left(\F^{\op m}\right)=c\left(\F\right)^{-m}$.
\end{proof}

\begin{prop}\label{prop:FormulaCasoK1}
Set
\begin{equation}\label{eq:PolinomioCasoK1}
P_{n,r}\left(m\right)=\sum_{i=0}^{n-1}\left(-1\right)^{i}\left(2r\right)^{n-1-i}\binom ni
\sum_{\substack{u,v\ge0\\u+v=n-1-i}}\binom{m+u-1}{u}\binom mv\ \in\QQ\left[m\right] .
\end{equation}
Then $\deg P_{n,r}=n-1$ and $\ged\left(X_{m}\right)=P_{n,r}\left(m\right)$ for every $m\ge n$.
\end{prop}

\begin{proof}
Substituting Lemma \ref{lem:FibratiCasoK1} into \eqref{eq:gEDFormulaBuona} and using $i+u+v=b=n-1$ in every summand gives \eqref{eq:PolinomioCasoK1}; the degree is $n-1$ by Theorem \ref{thm:PolinomialitaStabileInM}. Since $k+b=n$, the equality $\ged\left(X_{m}\right)=P_{n,r}\left(m\right)$ for $m\ge n$ is Theorem \ref{thm:PolinomialitaStabileInM}.
\end{proof}

\begin{example}\label{ex:FormuleEsplicite}
For $n=2,3,4$ formula \eqref{eq:PolinomioCasoK1} gives
\[
P_{2,r}\left(m\right)=4rm-2,\qquad P_{3,r}\left(m\right)=8r^{2}m^{2}-12rm+3 ,
\]
\[
P_{4,r}\left(m\right)=\tfrac{1}{3}\left(32r^{3}m^{3}-96r^{2}m^{2}
+\left(16r^{3}+72r\right)m-12\right).
\]
\end{example}

Since $X_{m}$ is smooth, \eqref{eq:FormulaChernLiscio} computes $\ged\left(X_{m}\right)$, and we may compare it with $P_{n,r}\left(m\right)$. In the table below each entry is $\ged\left(X_{m}\right)$; when $P_{n,r}\left(m\right)$ differs from it, the two are separated by a slash, $\ged\left(X_{m}\right)$ first.

\begin{center}
\begin{tabular}{cc|ccccc}
$n$ & $r$ & $m=1$ & $m=2$ & $m=3$ & $m=4$ & $m=5$\\
\hline
$2$ & $2$ & $4\,/\,6$ & $14$ & $22$ & $30$ & $38$\\
$3$ & $2$ & $13\,/\,11$ & $79\,/\,83$ & $219$ & $419$ & $683$\\
$3$ & $3$ & $39$ & $210\,/\,219$ & $543$ & $1011$ & $1623$\\
$4$ & $3$ & $203\,/\,212$ & $1562\,/\,1580$ & $5801\,/\,5828$ & $14684$ & $29876$\\
\end{tabular}
\end{center}

\begin{rem}\label{rem:SogliaOttimale}
The two computations agree for $m\ge n$, as they must, and disagree for $m<n$, except in the case $\left(n,r,m\right)=\left(3,3,1\right)$, where they happen to coincide. Hence the bound $m\ge k+b$ of Theorem \ref{thm:PolinomialitaStabileInM} cannot be lowered uniformly, although it may fail to be optimal for particular triples, see Remark \ref{rem:SogliaNonSempreOttimale}.
\end{rem}

\begin{rem}\label{rem:ConfrontoKMQS}
For $r=1$ the variety $X_{m}$ is the variety of $m\times n$ matrices of rank one, whose generic $\ED$-degrees are tabulated in \cite[Table 2]{KMQS25} for $n-1\le m-1\le10$: formula \eqref{eq:PolinomioCasoK1} reproduces all of them, the first rows being $P_{2,1}\left(m\right)=4m-2$, $P_{3,1}\left(m\right)=8m^{2}-12m+3$ and $P_{4,1}\left(m\right)=\frac{1}{3}\left(32m^{3}-96m^{2}+88m-12\right)$, and the range $n\le m$ in which they are tabulated is exactly our stable range. Note that $P_{n,r}$ is polynomial in $\left(m,r\right)$ jointly of degree $n-1$ in each, as predicted by Theorem \ref{thm:PolinomialitaCongiunta}, with leading coefficient $4^{n-1}r^{n-1}/\left(n-1\right)!$: this agrees with Theorem \ref{thm:PolinomialitaStabileInM}, the degree of $v_{r}\PP\left(V\d\right)$ in $\PP\left(S_{r}\right)$ being $r^{n-1}$.
\end{rem}

\subsection{The case $k=2$}\label{sec:esempiBinario}%~\\
Throughout this subsection we fix $k=2$ and $r\ge2$, so that Proposition \ref{prop:BaseFissa} applies and $\Sigma\simeq\Hilb_{2}\PP\left(V\d\right)$, which is smooth, as recalled in Section \ref{sec:baseFissa}. In particular one can take $B=\Sigma$ and $\T=T_{\Sigma}$: neither the Nash blow-up nor a resolution is needed, and \eqref{eq:gEDFormulaBuona} becomes an integral of Chern classes over a smooth variety. For $n=2$ that variety is $\PP^{2}$ and the integral can be evaluated explicitly. In the case of arbitrary $n$ we content ourselves with the leading term of $P_{n,2,r}$, which by Theorem \ref{thm:PolinomialitaStabileInM} only requires the degree of $\Sigma_{2,r}^{n}$. We point out that the case $n=k=r=2$ is also where our results meet the determinantal computation of \cite[Theorem 6.5]{KLW24}, discussed in Remark \ref{rem:ConfrontoKLW}.

\begin{lem}\label{lem:FibratoCasoBinario}
Let $n=2$. Then $\Sigma\simeq\Hilb_{2}\PP^{1}\simeq\PP^{2}$ is smooth, so one may take $B=\Sigma$ and $\T=T_{\PP^{2}}$. Let $\gamma$ denote the hyperplane class, then the
tautological subbundle $\EE$ of Section \ref{sec:costruzioneCm} satisfies
\[
c\left(\EE\right)=\left(1+\gamma\right)^{-\left(r-1\right)} .
\]
\end{lem}

\begin{proof}
By Proposition \ref{prop:BaseFissa}, $\Sigma \simeq \Hilb_{2}\PP^{1} \simeq \PP\left(S_{2}\d\right) \simeq \PP^{2}$ identifies with the space of binary quadrics; it is smooth, so its Nash blow-up is the identity. The universal family is $\mathscr{Z}=\left\{ \left(x,q\right)\,\vert\,q\left(x\right)=0\right\} \su\PP^{1}\times\PP^{2}$ of bidegree $\left(2,1\right)$. Twisting $0\to\OO\left(-2,-1\right)\to\OO\to\OO_{\mathscr{Z}}\to0$ by $\OO\left(r,0\right)$ and pushing forward to $\PP^{2}$ gives, for $r\ge2$,
\[
0\to\OO\left(-1\right)^{\op\left(r-1\right)}\to\OO^{\op\left(r+1\right)}\to\mathcal{F}_{r}\to0,
\]
whence $c\left(\mathcal{F}_{r}\right)=\left(1-\gamma\right)^{-\left(r-1\right)}$ and $c\left(\EE\right)=c\left(\mathcal{F}_{r}\d\right)=\left(1+\gamma\right)^{-\left(r-1\right)}$.
\end{proof}

\begin{prop}\label{prop:FormulaCasoBinario}
Let $n=2$ and $m\ge4$. Then
\[
\ged\left(X_{2,r}^{2,m}\right) = 8\left(r-1\right)^{2}m^{2}-12\left(r-1\right)m+3 .
\]
\end{prop}

\begin{proof}
Here $b=2$ and $\int_{\PP^{2}}\gamma^{2}=1$. By Lemma \ref{lem:FibratoCasoBinario} and Whitney's formula $c\left(\Q\d\right)=\left(1-\gamma\right)^{r-1}$, so that
\[
s_{u}\left(\EE^{\op m}\right)=\binom{m\left( r-1\right)}{u}\gamma^{u},\qquad
s_{v}\left(\left(\Q\d\right)^{\op m}\right)=\binom{m\left( r-1\right)+v-1}v\gamma^{v},
\]
while $c\left(\T\d\right)=\left(1-\gamma\right)^{3}$. Substituting in \eqref{eq:gEDFormulaBuona}, the three contributions are
\[
i=0\to 4\Big[\tbinom {m\left( r-1\right)}{2}+m^2\left( r-1\right)^{2}+\tbinom{m\left( r-1\right)+1}2\Big]=8m^2\left( r-1\right)^{2}
\]
\[
i=1\to -2\cdot3\cdot2m\left( r-1\right)=-12m\left( r-1\right),\qquad i=2\to 3 .
\]
\end{proof}

\begin{rem}\label{rem:CasoBDue}
Proposition \ref{prop:FormulaCasoBinario} and the case $n=3$ of Example \ref{ex:FormuleEsplicite} both give the value $8t^{2}-12t+3$, for $t=m\left(r-1\right)$ and $t=mr$ respectively. This is not a coincidence: whenever $b=2$ and $B\simeq\PP^{2}$ with $\T=T_{\PP^{2}}$, Whitney's formula forces $c_{2}\left(\EE\right)+c_{2}\left(\Q\d\right)=c_{1}\left(\EE\right)^{2}$ and the second Chern classes cancel from \eqref{eq:gEDFormulaBuona}, so that the answer depends only on $c_{1}\left(\EE\right)$.
\end{rem}

\begin{prop}\label{prop:SogliaRaffinataKLW}
For $n=k=r=2$ the conclusion of Proposition \ref{prop:FormulaCasoBinario} holds for every $m\ge3$, and fails for $m=2$.
\end{prop}

\begin{proof}
For $m=3$ one has $m<k+b=4$, so Proposition \ref{prop:ProprietaCm} and Theorem \ref{thm:PolinomialitaStabileInM} do not apply as stated and we re-examine their proofs. Besides the standing assumption $m\ge k$, the hypothesis $m\ge k+b$ is used there in exactly three places: to guarantee $m\left(N-k\right)>b$, which is the hypothesis of Lemma \ref{lem:ConormaleSuXcirc} and gives $\rk\K_{m}>0$; to guarantee, through the bound $\mathrm{codim}_{\P_{m}}\tilde{\D}_{m}=m-k+1>b$, that no irreducible component of $\C_{m}$ is contained in $\tilde{\D}_{m}$; and to guarantee \eqref{eq:CondizioneCombinatoria} for every $g\le b$ in Lemma \ref{lem:Coefficiente2g}. We verify the three conditions directly for $n=k=r=2$ and $m\ge3$, so that both proofs go through verbatim and, with them, Proposition \ref{prop:FormulaCasoBinario}.

Here $N=3$ and $b=2$. The first condition reads $m\left(3-2\right)>2$, that is $m\ge3$; the third, for $g\le2$, reads $2\le2m-1$ and $2\le m-1$, again $m\ge3$. As for the second, we bound $\dim\left(\C_{m}\cap\tilde{\D}_{m}\right)$ directly. A point of $\tilde{\D}_{m}$ has $\rk T^{\flat}\le1$ and $T\neq0$, hence $\rk T^{\flat}=1$; fibrewise over $B$ the rank one locus of $\PP\left(\Hom\left(W_{m}\d,\EE_{z}\right)\right)$ is the Segre variety $\PP^{m-1}\times\PP^{k-1}$, of dimension $m$, so that the locus $\left\{ \rk T^{\flat}=1\right\} \su\H_{m}$ has dimension $b+m=m+2$. Fix such a point and let $e$ span $\mathrm{Im}\left(T^{\flat}\right)$: the map $\lambda_{z,T}$ factors through $\mathrm{Im}\left(T^{\flat}\right)\ot\Q\d_{z}$, of dimension $N-k=1$, so $\rk\lambda_{z,T}\le1$; and $\lambda_{z,T}=0$ would force $\xi\left(e\right)=0$ for every $\xi\in\T_{z}=T_{\EE_{z}}\Sigma$, impossible since $\Sigma=\Gr_{2}\left(S_{2}\right)$ and hence $\T_{z}=\Hom\left(\EE_{z},\Q_{z}\right)$. Therefore $\rk\lambda_{z,T}=1$, the fibre of $\C_{m}$ over $\left(z,\left[T\right]\right)$ is a linear subspace of $\PP\left(W_{m}\d\ot\Q\d_{z}\right)\simeq\PP^{m-1}$ of dimension $m-2$, and 
\[
\dim\left(\C_{m}\cap\tilde{\D}_{m}\right)=\left(m+2\right)+\left(m-2\right)=2m ,
\]
which is smaller than $\dim\C_{m}=mN-2=3m-2$ precisely when $m\ge3$.

For $m=2$ all three conditions fail, and the conclusion fails as well: every $2\times3$ matrix has rank at most $2$, so $X_{2,2}^{2,2}=\PP^{5}$, which is excluded from Definition \ref{def:DefinizioneGradiPolariConormale}, its conormal variety being empty, while Proposition \ref{prop:FormulaCasoBinario} would give $11$. The three failures are consistent with one another: the first condition fails precisely because the architecture is  filling, and the count above gives $\dim\left(\C_{2}\cap\tilde{\D}_{2}\right)=4=\dim\C_{2}$, so that $\C_{2}$ acquires a component contained in $\tilde{\D}_{2}$ and $\left[\C_{2}\right]$ no longer pushes forward to $\left[\N_{X_{2}}\right]$.
\end{proof}

\begin{rem}\label{rem:SogliaNonSempreOttimale}
Proposition \ref{prop:SogliaRaffinataKLW} shows that the bound $m\ge k+b$ of Theorem \ref{thm:PolinomialitaStabileInM} is sufficient but not always optimal: for $n=k=r=2$ the true threshold is $m\ge3$ against $k+b=4$, and it is exactly $3$ because the three conditions of Remark \ref{rem:DueCondizioni} become binding at the same value of $m$. It cannot, however, be replaced by any smaller expression in $\left(n,k,r\right)$: for $k=1$ one has $k+b=n$, and in each of the cases displayed in the table of Section \ref{sec:esempiK1} one has $\ged\left(X_{m}\right)\neq P_{n,r}\left(m\right)$ at $m=n-1$, so that for those triples $k+b$ is exactly the threshold.
\end{rem}

\begin{rem}\label{rem:ConfrontoKLW}
The architecture  $\left(2,2,d_{2}\right)$ with $r=2$ of \cite{KLW24}, that is the network $\CC^{2}\to\CC^{2}\to\CC^{d_{2}}$, corresponds to $\left(n,k,r\right)=\left(2,2,2\right)$ and $m=d_{2}$, so that $N=3$ and $X_{2,2}^{2,m}\su\PP^{3m-1}$ is by \cite[Proposition 4.1]{KLW24} the determinantal variety of $m\times3$ matrices of rank at most $2$. At $r=2$ Proposition \ref{prop:FormulaCasoBinario} gives $8m^{2}-12m+3$, that is \cite[Theorem 6.5]{KLW24}, where the same value is obtained through the Chern-Mather classes of determinantal varieties and the formula of \cite[Proposition 5.5]{Zha18}: the determinantal case is thus the slice $r=2$ of a one-parameter family, and our derivation does not use the determinantal description.

However \cite[Theorem 6.5]{KLW24} states the equality for every $d_{2}\ge2$, that is for every $m\ge2$ in our notation, whereas by Proposition \ref{prop:SogliaRaffinataKLW} it can only hold for $m\ge3$: for $m=2$ the architecture is filling and $X_{2,2}^{2,2}=\PP^{5}$, by \cite[Proposition 4.1]{KLW24}.

The correct value there is $\ED_q\left(X_{2,2}^{2,2}\right)=1$ for every $q$: the affine cone is the whole of $\M_{2}$, and the squared distance from a general $u$ has $u$ itself as its only critical point.  Note that the $\ged$ of \eqref{eq:DefinizioneGEDGradiPolari} is not defined here, the conormal variety being empty. The discrepancy can be traced back to the source of that formula: in the notation of \cite{Zha18} the variety of $m\times3$ matrices of rank at most $2$ is $\tau_{m,3,1}$, the locus of matrices whose kernel has dimension at least one, and the Chern-Mather class on which \cite[Proposition 5.5]{Zha18} rests is computed in \cite[Theorem 4.3]{Zha18} under the standing hypothesis $n\le m$, which here reads $m\ge3$: exactly the threshold of Proposition \ref{prop:SogliaRaffinataKLW}. We stress that the consequence drawn in \cite[Theorem 6.5]{KLW24}, namely that the learning degree is at most $8m^{2}-12m+3$, remains correct for $m=2$ as well.
\end{rem}

\begin{rem}\label{rem:DualitaKLW}
The polynomial $8m^{2}-12m+3$ is also $P_{3,1}\left(m\right)$ of Example \ref{ex:FormuleEsplicite}, and this is not a coincidence:  $X_{1,1}^{3,m}=\Seg\left(\PP\left(W_{m}\right)\times\PP\left(V\d\right)\right)$ is the variety of $m\times3$ matrices of rank one, whose projective dual is, for $m\ge3$, the variety of $m\times3$ matrices of rank at most $2$, that is $X_{2,2}^{2,m}$; and $\ged\left(X\right)=\ged\left(X\d\right)$ by Proposition \ref{prop:ProprietaGed}$\left(1\right)$. Combined with Remark \ref{rem:ConfrontoKMQS}, this shows that the value also appears, in the dual guise, in \cite[Table 2]{KMQS25}.
\end{rem}

By Theorem \ref{thm:PolinomialitaCongiunta} the leading coefficient of $P_{n,2,r}$ is governed by $\deg\left(\Sigma_{2,r}^{n}\right)$, which we now compute. By Proposition \ref{prop:BaseFissa} the base is $\Hilb_{2}\PP\left(V\d\right)$ and $c_{1}\left(\Q\right)=c_{1}\left(\EE_{r}\right)$, where $\EE_{r}=p_{\ast}\xi\inv\OO\left(r\right)$ is the tautological rank two bundle attached to $\OO\left(r\right)$: the degree we are after is thus a tautological Chern number on a Hilbert scheme of two points. In \cite[Lemma 2.1]{Sta26} the author computes
\begin{equation}\label{eq:FormulaIntegrale}
\int_{S^{\left[2\right]}}c_{1}\left(\mathcal{L}^{\left[2\right]}\right)^{2d}=
\frac{1}{2}\binom{2d}{d} \Big(\int_{S}\left( c_1\left(\mathcal L\right)\right)^{d}\Big)^{2} -2^{d-1}\sum_{j=0}^{d}\binom{2d}{d+j}\frac{1}{2^{j}} \int_{S}\left(c_1\left(\mathcal L\right)\right)^{d-j}s_{j}\left(S\right) 
\end{equation}
for every smooth projective $S$ of dimension $d$ and every invertible sheaf $\mathcal{L}$ on it, in terms of $\int_{S}c_{1}\left(\mathcal{L}\right)^{d}$ and of the Segre classes of $S$. Our case is $S=\PP^{n-1}$ and $\mathcal{L}=\OO\left(r\right)$, and all that is left is to evaluate the resulting expression.

\begin{lem}\label{lem:GradoSigma}
For $n\ge2$ and $r\ge2$,
\[
\deg\left(\Sigma_{2,r}^{n}\right)=\binom{2n-3}{n-1}\Big[r^{2n-2}-\left(2r-1\right)^{n-1}\Big] .
\]
\end{lem}

\begin{proof}
Let $P=\PP\left(V\d\right)$, identify $\Sigma$ with $\Hilb_{2}P$ in view of Proposition \ref{prop:BaseFissa} and put $d=n-1$. Then $\EE_{r}$ identifies with the tautological rank two bundle $\OO_{P}\left(r\right)^{\left[2\right]}$ and $c_{1} \left(\Q\right) = c_{1} \left(\EE_{r}\right)$, so that
\[
\deg\left(\Sigma_{2,r}^{n}\right)
=\int_{\Hilb_{2}P}c_{1}\left(\OO_{P}\left(r\right)^{\left[2\right]}\right)^{2d}
\]
and \eqref{eq:FormulaIntegrale} applies: writing $h$ for the hyperplane class of $P$,
\[
\deg\left(\Sigma_{2,r}^{n}\right)=\frac{1}{2}\binom{2d}{d} \Big(\int_{P}\left(rh\right)^{d}\Big)^{2} -2^{d-1}\sum_{m=0}^{d}\binom{2d}{d+m}\frac{1}{2^{m}} \int_{P}\left(rh\right)^{d-m}s_{m}\left(T_{P}\right) .
\]
The first term is $\frac{1}{2}\binom{2d}{d}r^{2d}$. As for the second, the Euler sequence gives $s_{m}\left(T_{P}\right)=\left(-1\right)^{m}\binom{d+m}{m}h^{m}$, whence $\int_{P}\left(rh\right)^{d-m}s_{m}\left(T_{P}\right)=\left(-1\right)^{m}\binom{d+m}{m}r^{d-m}$. Using $\binom{2d}{d+m}\binom{d+m}{m}=\binom{2d}{d}\binom{d}{m}$ the sum becomes
\[
2^{d-1}\binom{2d}{d}\sum_{m=0}^{d}\binom{d}{m}\left(-\frac{1}{2}\right)^{m}r^{d-m}
=2^{d-1}\binom{2d}{d}\left(r-\frac{1}{2}\right)^{d}
=\frac{1}{2}\binom{2d}{d}\left(2r-1\right)^{d}.
\]
Therefore $\deg\left(\Sigma_{2,r}^{n}\right)=\frac{1}{2}\binom{2d}{d} \left[r^{2d}-\left(2r-1\right)^{d}\right]$, and the assertion follows from $\frac{1}{2}\binom{2d}{d}=\binom{2d-1}{d}$ with $d=n-1$.
\end{proof}

\begin{cor}\label{cor:TermineDominanteKDue}
For $n\ge2$, $r\ge2$ and $m\ge2n$,
\[
\ged\left(X_{2,r}^{n,m}\right)=\frac{4^{2n-2}}{\left(2n-2\right)!}\binom{2n-3}{n-1}
\Big[r^{2n-2}-\left(2r-1\right)^{n-1}\Big]m^{2n-2}+O\left(m^{2n-3}\right).
\]
\end{cor}

\begin{proof}
Lemma \ref{lem:GradoSigma} and Theorem \ref{thm:PolinomialitaStabileInM}, with $b=2\left(n-1\right)$.
\end{proof}

\newpage
\bibliographystyle{amsplain}
\bibliography{References}
\end{document}